\documentclass[onecolumn]{IEEEtran}
\usepackage{cite}
\usepackage{amsmath,amssymb,amsfonts}
\usepackage{algorithmic}
\usepackage{graphicx}
\usepackage{textcomp}
\usepackage{mathrsfs}

\usepackage{epsfig} 
\usepackage{comment}
\usepackage{bm}
\usepackage{color}
\usepackage{algorithm}
\usepackage{stmaryrd}
\newtheorem{problem}{Problem}
\newtheorem{proposition}{Proposition}

\newtheorem{lemma}{Lemma}
\newtheorem{theorem}{Theorem}

\newtheorem{assumption}{Assumption}
\newtheorem{rmk}{Remark}

\newcommand{\argmin}{\mathop{\rm arg~min~}\limits}
\newcommand{\esssup}{\mathop{\rm ess~sup}\limits}

\newcommand{\blue}[1]{{\color{blue}{#1}}}
\newcommand{\red}[1]{{\color{red}{#1}}}

\def\BibTeX{{\rm B\kern-.05em{\sc i\kern-.025em b}\kern-.08em
    T\kern-.1667em\lower.7ex\hbox{E}\kern-.125emX}}
\begin{document}
\title{Sparse optimal control for infinite-dimensional linear systems with applications to graphon control}
\author{Takuya Ikeda and Masaaki Nagahara
\thanks{This work was supported in part by the Japan Society for the Promotion of Science (JSPS) KAKENHI under Grant JP24K17300, 23K26130, 22KK0155, 22H00512, and 24K21314.
This work is also supported by Japan Science and Technology Agency (JST) as part of Adopting Sustainable Partnerships for Innovative Research Ecosystem (ASPIRE), Grant Number JPMJAP2402.}
\thanks{T. Ikeda is with the Faculty of Environmental Engineering, University of Kitakyushu, Fukuoka, 808-0135, Japan
(e-mail: t-ikeda@kitakyu-u.ac.jp). }
\thanks{M. Nagahara is with the Graduate School of Advanced Science and Engineering, Hiroshima University, 739-8521, Japan
(e-mail: nagam@hiroshima-u.ac.jp). }
}

\maketitle

\begin{abstract}
Large-scale networked systems typically operate under resource constraints, and it is also difficult to exactly obtain the network structure between nodes.
To address these issues, this paper investigates a sparse optimal control for infinite-dimensional linear systems 
and its application to networked systems where the network structure is represented by a limit function called a graphon that captures the overall connection pattern.
The contributions of this paper are twofold:
(i) To reduce computational complexity, we derive a sufficient condition under which the sparse optimal control can be obtained by solving its corresponding $L^1$ optimization problem. 
Furthermore, we introduce a class of non-convex optimal control problems such that the optimal solution always coincides with a sparse optimal control, provided that the non-convex problems admit optimal solutions.
(ii) We show that the sparse optimal control for large-scale finite-dimensional networked systems can be approximated by that of the corresponding limit graphon system, provided that the underlying graph is close to the limit graphon in the cut-norm topology.
The effectiveness of the proposed approach is illustrated through numerical examples.\\


\end{abstract}


\section{Introduction}

The study of large-scale networked systems has become increasingly important across a wide range of research fields, including engineering, economics, biology, and the social sciences.
These systems are often subject to practical constraints such as limited actuation, communication bandwidth, and energy availability, which make it inefficient or even infeasible to control or monitor all nodes at all times.
Accordingly, there is a growing need for the development of \emph{resource-aware control} that can achieve high performance while utilizing minimal control resources.

\emph{Sparse optimal control} has been proposed as a promising approach to this challenge, which aims to minimize the number of active control inputs while achieving the desired control performance~\cite{NagQueNes16}. 
For example, sparse control techniques are employed in actuator placement problems to identify a small subset of nodes 
that receive control inputs to effectively guide the overall system~\cite{ManKutBru21, IkeKas22}. 
Furthermore, sparse control naturally results in extended periods during which actuators remain inactive, thereby contributing to reductions in fuel consumption, power usage, and communication burden~\cite{Nag2020, Jos22}.

The natural penalty function for promoting sparsity is the $L^0$ norm, which measures the support length of a function.
However, due to the non-convex and discontinuous nature of the $L^0$ norm, a convex relaxation is typically employed, in which the $L^1$ norm is used as a convex surrogate of the $L^0$ norm.
It has been shown that, under certain conditions on the system model, the original sparse optimal control can be exactly obtained by solving its corresponding $L^1$ optimization problem~\cite{NagQueNes16, ITKKTAC18}. 
On the other hand, it is also known that the $L^1$ optimization does not always yield sparse solutions~\cite{AthFal}. 
To address this limitation, recent studies have investigated non-convex optimization problems that promote the sparsity under less restrictive conditions~\cite{HayIkeNag24, Ike25}.

The aforementioned studies focus on sparse optimal controls for given finite-dimensional systems.
However, to effectively extend the methods to large-scale networked systems, it is necessary to address some additional challenges inherent to the large-scale nature.
For example, it is difficult or too costly to exactly obtain the network structure between nodes.
Furthermore, control designs that depend on the exact system size may require re-modeling and re-computation of optimal controls even for slight changes in size, which can lead to undesirable consumption of computational resources and energy. 
For the reasons, this paper investigates a sparse optimal control for infinite-dimensional systems and 
provides a scalable approximation method for sparse control in large-scale finite-dimensional networked systems based only on the overall connection pattern, while avoiding detailed description of individual interactions.

For the representation of network structures, this paper leverages \emph{graph limit theory} introduced and developed by~\cite{lovSze06, borChaLov08, borChaLov12, lov12}. 
This theory provides a rigorous framework for analyzing the limits of graph sequences when the number of nodes tends to infinity, and states that the limits can be represented by {\em graphons}.
A graphon is a bounded symmetric measurable function of two variables, and can be considered as a weighted graph whose underlying node set has continuum cardinality.
Its applications can be found in various fields, including 
mean field games~\cite{CaiHua19, BayCha23}, 
signal processing~\cite{MorLeu17, RuiCha21}, 
neural networks~\cite{RuiCha20, Levi23},
and epidemics~\cite{VizFra20, KelHor22}.
In a relevant study~\cite{GaoCai20}, the $L^2$ optimal control problem with unbounded control inputs for systems defined by graphons is considered, and an approximate control strategy for finite-dimensional systems is proposed based on step function approximations and a closed-form expression of the optimal control.
However, to the best of our knowledge, no existing study has addressed the sparsity of control in graphon-based formulations. 

Our contributions are as follows.
We first analyze a sparse optimal control problem for infinite-dimensional linear systems, which include networked systems defined by graphons.
We derive a sufficient condition under which the corresponding $L^1$ optimal control coincides with a sparse optimal control.
We also introduce a class of non-convex optimal control problems such that the optimal solution always coincides with a sparse optimal control.
Furthermore, we show that the sparse optimal control for large-scale finite-dimensional networked systems can be approximated by that of the corresponding limit graphon system, provided that the underlying graph is close to the limit graphon in the cut-norm topology, which is commonly used to measure structural similarity between large graphs.

The remainder of this paper is organized as follows.  
Section~\ref{sec:math} provides the mathematical preliminaries.  
Section~\ref{sec:graphon_system} formulates the sparse optimal control problem and analyzes the corresponding $L^1$ optimal controls and non-convex optimal controls.
Section~\ref{sec:approximation} shows the approximation result of the sparse control for finite-dimensional networked systems.  
Section~\ref{sec:example} illustrates the proposed method through numerical examples.  
Finally, Section~\ref{sec:conclusion} concludes the paper.

\section{Mathematical preliminaries}\label{sec:math}

This section reviews the notation that will be used throughout the paper.
The set of all integers is denoted by $\mathbb{Z}$,
the set of all positive integers is denoted by $\mathbb{N}$,
the set of all real numbers is denoted by $\mathbb{R}$,
the set of all non-negative numbers is denoted by $\mathbb{R}_{\geq0}$,
and the set $\{1, 2, \dots, n\}$ for some $n\in\mathbb{N}$ is denoted by $\left\llbracket n \right\rrbracket$.
For any $m\in\mathbb{N}$ and $\Omega\subset\mathbb{R}$,
$a=\left[a_1, a_2, \dots, a_m\right]^{\top}\in \Omega^{m}$ signifies that $a_i\in \Omega$ holds for all $i \in \left\llbracket m \right\rrbracket$.
The indicator function of a set $\Omega$ is denoted by ${\bf 1}_{\Omega}$.
The $\ell^p$ norm of $a\in\mathbb{R}^m$ for $p\in(0, \infty)$ is defined by
$\|a\|_{\ell^p} = \left( \sum_{i=1}^{m} |a_i|^p\right)^{\frac{1}{p}}$.
The Lebesgue measure on ${\mathbb{R}}$ is denoted by $\mu$. 
For a measurable function $u(t)=\left[u_1(t), u_2(t), \dots, u_m(t)\right]^{\top}\in{\mathbb{R}}^m$ over a measurable set $E$, the $L^p$ norm is defined by
\begin{align*}
	&\|u\|_{0} = \sum_{j=1}^{m}\mu(\{t\in E: u_j(t)\neq 0\}),\\
	&\|u\|_{p} = \left(\sum_{j=1}^{m}\int_{E} |u_j(t)|^p dt\right)^{\frac{1}{p}} \quad p\in(0,\infty),\\
	&\|u\|_{\infty} = \max_{1 \leq j \leq m} \esssup_{t \in E} |u_j(t)|.
\end{align*}
The set of all functions that have a finite $L^p$ norm on a measurable set $E\subset\mathbb{R}$ is denoted by $L^p_{E}$.
A property is said to hold almost everywhere (a.e.) if it holds everywhere except on some null set.
The closed ball with center $x \in L^2_{[0,1]}$ and radius $r>0$ is denoted by ${\mathrm{Ball}}(x, r)$, 
i.e., $\mathrm{Ball}(x, r)=\left\{y \in L^{2}_{[0,1]}: \left\|y - x\right\|_{2} \leq r\right\}$.
The dual space of $L^2_{[0,1]}$ is denoted by ${L^2_{[0,1]}}^\ast$, and the value of $f \in {L^2_{[0,1]}}^\ast$ at $x \in L^2_{[0,1]}$ is denoted by $\left< f, x \right> \in \mathbb{R}$.
An inner product on $L^2_{[0,1]}$ is denoted by $(\cdot, \cdot)$.
For a linear operator $A: L^2_{[0,1]} \to L^2_{[0,1]}$, its norm is denoted by $\left\|A\right\|_{\mathrm{op}}$, i.e., 
$\left\|A\right\|_{\mathrm{op}} = \sup_{ x \in L^2_{[0,1]} \backslash \{0\} } \frac{\left\|A x\right\|_{2}}{\left\|x\right\|_2}$.
The set of all bounded linear operators from $L^2_{[0,1]}$ to $L^2_{[0,1]}$ is denoted by $\mathcal{L}(L^2_{[0,1]})$.
The dual operator of $A \in \mathcal{L}(L^2_{[0,1]})$ is denoted by $A^\ast$, 
i.e., $\left<A^\ast f, x \right> = \left<f, Ax \right>$ for all $f \in {L^2_{[0,1]}}^\ast$ and $x \in L^2_{[0,1]}$.
The adjoint operator of $A\in \mathcal{L}(L^2_{[0,1]})$ is denoted by $A' $,
i.e., $\left(Ax, y\right) = \left(x, A'y\right)$ for all $x, y \in L^2_{[0,1]}$.
An operator-valued function $S(t)$ from $\mathbb{R}_{\geq0}$ to $\mathcal{L}(L^2_{[0,1]})$ is a {\em semigroup} if it satisfies
\begin{align*}
	&S(t+s) = S(t)S(s) \quad \mbox{for} \quad \forall t, \forall s \geq 0, \\
	&S(0) = I,
\end{align*}
where $I\in\mathcal{L}(L^2_{[0, 1]})$ denotes the identity operator.
A semigroup $S(t)$ is a {\em strongly continuous semigroup} if it satisfies
\begin{align*} 
	\lim_{t \to 0+} \left\|S(t)x - x \right\|_{2} = 0 \quad \mbox{for} \quad \forall x \in L^2_{[0,1]},
\end{align*}
and is a {\em uniformly continuous semigroup} if it satisfies
\begin{align*} 
	\lim_{t \to 0+} \left\|S(t) - I \right\|_{\mathrm{op}} = 0.
\end{align*}
(If $S(t)$ is a uniformly continuous semigroup, then it is clearly a strongly continuous semigroup.)
The {\em infinitesimal generator} $A$ of a semigroup $S(t)$ is defined by
\begin{align*}
	&Ax = \lim_{t\to0+} \frac{S(t)x - x}{t} \quad \mbox{for} \quad \forall x \in D(A),\\
	&D(A) = \left\{x \in L^2_{[0,1]}: \lim_{t\to0+} \frac{S(t)x - x}{t} \mbox{~exists}\right\}.
\end{align*}
For any $A\in\mathcal{L}(L^2_{[0,1]})$, $e^{At}$ is defined by
$e^{At} = \sum_{n=0}^{\infty} \frac{(At)^{n}}{n!}$.
It is known that 
$e^{At}$ is a uniformly continuous semigroup 
and $A$ is the infinitesimal generator of $e^{At}$ with $D(A) = L^2_{[0,1]}$ 
(see~\cite[Example 2.1.3 and Example 2.1.12]{CurZwa20}).

%

\section{Sparse optimal control for infinite-dimensional systems} \label{sec:graphon_system}

\subsection{Problem formulation}\label{subsec:formulation_1}

In this section, we consider an infinite-dimensional linear system given by
\begin{equation}\label{eq:system}
	\dot{x}(t)=Ax(t)+Bu(t), \quad 0 \leq t \leq T,
\end{equation}
where 
$x(t)\in L^2_{[0,1]}$ is the state at time $t$ with an initial state $x_{0} \in L^{2}_{[0,1]}$,
$u(t) =\left[u_{1}(t), u_{2}(t), \dots, u_{m}(t)\right]^\top \in {\mathbb{R}}^m$ is the control input at time $t$,
$T>0$ is the final time of control,
and $A\in\mathcal{L}(L^2_{[0,1]})$ is a bounded linear operator.
The bounded linear operator $B: \mathbb{R}^m \to L^2_{[0,1]}$ is defined by $Bu = \sum_{j=1}^{m} b_j u_{j}$ with $b_1, b_2, \dots, b_m \in L^2_{[0,1]}$, and we denote by $\mathscr{B}$ the set of all such operators.
The system~\eqref{eq:system} is simply denoted by $(A; B)$ or $(A; B; x_0)$ throughout the paper.
For the system~\eqref{eq:system}, a control input $u$ is said to be {\em admissible} if it satisfies the constraint $\left\|u\right\|_{\infty} \leq 1$, and the set of all admissible controls is denoted by 
\[
	\mathcal{U} = \left\{u\in L_{[0, T]}^{\infty}: \left\|u\right\|_{\infty} \leq 1\right\}.
\]
In this setting, the state is the mild solution of~\eqref{eq:system} given by
\begin{equation}\label{eq:mild solution}
	x(t) = e^{At} x_0 + \int_{0}^{t} e^{A(t-s)} B u(s) ds, \quad 0 \leq t \leq T,
\end{equation}
where $e^{At}$ is the uniformly continuous semigroup with infinitesimal generator $A$ (see \cite[p. 190]{CurZwa20}).
Note that \eqref{eq:mild solution} is a well defined integral in the sense of Bochner (see \cite[Lemma A.5.10 and Example A.5.20]{CurZwa20}).

Here, the {\em sparse optimal control} refers to the optimal solutions to the following optimal control problem, where the parameter $\lambda>0$ serves as a trade-off factor, balancing the sparsity of the control input $u$ and the error between the terminal state $x(T)$ and the target state $x_{f}$:

\begin{problem}[sparse optimal control]\label{prob:main_limit}
For $A\in\mathcal{L}(L^2_{[0,1]})$, $B \in \mathscr{B}$, $x_{0}, x_{f}\in L^{2}_{[0,1]}$, $T>0$, and $\lambda>0$, 
find a control input $u$ on $[0, T]$ that solves the following:
\begin{equation*}
\begin{aligned}
  & \underset{u}{\text{minimize}}
  & & \left\|u\right\|_{0} + \lambda \left\|x(T) - x_{f}\right\|_{2}^2\\
  & \text{subject to}
  & & \dot{x}(t) = Ax(t) + Bu(t), \quad x(0) = x_{0}, \quad u\in\mathcal{U}.
\end{aligned}
\end{equation*}
\end{problem}

Owing to the $L^0$ norm in the cost function, it is computationally difficult to solve Problem~\ref{prob:main_limit} as it is.
To circumvent this issue, this paper first considers a convex relaxation method that replaces the $L^0$ norm with the $L^1$ norm.
The problem is formulated as follows, where the optimal solution is called {\em the $L^1$ optimal control}:

\begin{problem}[$L^1$ optimal control]\label{prob:convex_limit}
For $A\in\mathcal{L}(L^2_{[0,1]})$, $B \in \mathscr{B}$, $x_{0}, x_{f}\in L^{2}_{[0,1]}$, $T>0$, and $\lambda>0$, 
find a control input $u$ on $[0, T]$ that solves the following:
\begin{equation*}
\begin{aligned}
  & \underset{u}{\text{minimize}}
  & & \left\|u\right\|_{1} + \lambda \left\|x(T) - x_{f}\right\|_{2}^2\\
  & \text{subject to}
  & & \dot{x}(t) = Ax(t) + Bu(t), \quad x(0) = x_{0}, \quad u\in\mathcal{U}.
\end{aligned}
\end{equation*}
\end{problem}

In Section~\ref{subsec:analysis_1}, it is shown that the $L^1$ optimal control is exactly a sparse optimal control under a condition on the system $(A; B)$ (see Theorem~\ref{thm:L0-L1}).
However, in general, sparse optimal control cannot always be obtained through $L^1$ optimization~\cite{ITKKTAC18}.
Then, in Section~\ref{subsec:analysis_2}, we also present non-convex optimal control problems which always give a sparse optimal control under an assumption on the cost function.
Throughout the paper, the cost functions of Problems~\ref{prob:main_limit} and \ref{prob:convex_limit} are denoted by $J_0: \mathcal{U} \to \mathbb{R}$ and $J_1: \mathcal{U} \to \mathbb{R}$, i.e., 
\begin{align*}
	&J_0(u) = \left\|u\right\|_{0} + \lambda \left\|x(T) - x_{f}\right\|_{2}^2,\\
	&J_1(u) = \left\|u\right\|_{1} + \lambda \left\|x(T) - x_{f}\right\|_{2}^2.
\end{align*}

\subsection{Analysis on $L^1$ optimal control}
\label{subsec:analysis_1}

We first show a necessary condition of the $L^1$ optimal control.

\begin{theorem}\label{thm:max-cond-convex}
Let $\check{u}$ be any optimal solution to Problem~\ref{prob:convex_limit} and $ \check{x}$ be the corresponding state.
Then, we have
\begin{equation}\label{eq:necessary_L1}
	\check{u}(t) \in \argmin_{v\in[-1,1]^m} 
	\left\{\|v\|_{\ell^1} +  2\lambda \left( \check{x}(T) - x_{f}, e^{A(T-t)} Bv \right) \right\}
\end{equation}
almost everywhere.
\end{theorem}

\begin{IEEEproof}
Note that all the assumptions in~\cite[Theorem 6.6.2]{Fat99} are satisfied in our case, as mentioned in~\cite[Remark 6.6.8]{Fat99}.
Then, it follows from~\cite[Theorem 6.6.2]{Fat99} that
there exists $(z_0, z)\in \mathbb{R}_{\geq0} \times {L_{[0,1]}^2}^{\ast}$ satisfying following: 
\begin{align}
	&z\in N_{L^2_{[0,1]}}(\check{x}(T)),\notag\\
	&\check{u}(t) \in \argmin_{v\in[-1,1]^m} \left\{z_0\left\|v\right\|_{\ell^1} + \left<\check{z}(t), Bv \right> \right\} \mbox{~a.e.~} t\in[0,T],\label{eq:min_cond}
\end{align}
where 
$\check{z}(t)={e^{A(T-t)}}^\ast(z+z_0 \partial g_0(\check{x}_T))$,
${L_{[0,1]}^2}^{\ast}$ is the dual space of $L^2_{[0,1]}$,
${e^{A(T-t)}}^\ast$ is the dual operator of $e^{A(T-t)}$, 
$N_{L^2_{[0,1]}}(\check{x}(T))$ is the normal cone to ${L^2}_{[0,1]}$ at $\check{x}(T)$,
$g_0: L^2_{[0,1]}\to\mathbb{R}$ is the functional defined by
$g_0(x) = \lambda \left\|x-x_{f}\right\|_2^2$,
and $\partial g_0(\check{x}(T))$ is the Fr\'{e}chet differential of $g_0$ at $\check{x}(T)$.
Here, we have
$N_{L^2_{[0,1]}}(\check{x}(T))=\{0\}$
and $\left< \partial g_0(\check{x}(T)), y \right> = 2\lambda \left(\check{x}(T) - x_{f}, y\right)$ for any $y\in L^2_{[0,1]}$.
Hence, we have
\begin{align}\label{eq:min_cond_dual}
\begin{split}
	\left<\check{z}(t), Bv \right>
	&= z_0 \left< {e^{A(T-t)}}^{\ast} \partial g_0(\check{x}(T)), Bv \right>\\
	&= z_0 \left< \partial g_0(\check{x}(T)), e^{A(T-t)} Bv \right>\\
	&= 2\lambda z_0 \left( \check{x}(T) - x_{f}, e^{A(T-t)} Bv \right).
\end{split}
\end{align}

We next show $z_0 > 0$. 
For this, let us take any ${u}^{\prime}_{(n)} \in \mathcal{U}$ such that $\left\| {u}^{\prime}_{(n)} - \check{u} \right\|_{0} \to 0$ as $n\to\infty$,
and denote by ${x}^{\prime}_{(n)}$ the state corresponding to ${u}^{\prime}_{(n)}$.
Define
\begin{align*}
	&R({u}^{\prime}_{(n)}) = 
	\left\{ \int_{0}^{T} e^{A(T-t)} B \left( u(t) - {u}^{\prime}_{(n)}(t) \right) dt: u \in \mathcal{U} \right\}, \\
	&{w}^{\prime}_{(n)} = \frac{1}{T} \int_{0}^{T} e^{A(T-t)} B {u}^{\prime}_{(n)}(t) dt, \\
	&\check{w} = \frac{1}{T} \int_{0}^{T} e^{A(T-t)} B \check{u}(t) dt.
\end{align*}
Note that
\begin{align*}
	\left\| {w}^{\prime}_{(n)} - \check{w} \right\|_{2}
	&\leq \frac{1}{T} \int_{0}^{T} \left\| e^{A(T-t)} B \left( {u}^{\prime}_{(n)}(t) - \check{u}(t) \right) \right\|_{2} dt\\
	&\leq \frac{1}{T} e^{\|A\|_{\mathrm{op}}T} \int_{0}^{T} \left\| B \left( {u}^{\prime}_{(n)}(t) - \check{u}(t) \right)\right\|_{2} dt\\
	&\leq \frac{2}{T} e^{\|A\|_{\mathrm{op}}T} \left\|{u}^{\prime}_{(n)} - \check{u} \right\|_{0} \sum_{j=1}^{m} \|b_j\|_{2},
\end{align*}
where the first relation follows from Minkowski's integral inequality.
Hence, 
for any $\varepsilon>0$, there exists $n_{\varepsilon}\in\mathbb{N}$ such that 
if $n\geq n_{\varepsilon}$, then we have $\left\| {w}^{\prime}_{(n)} - \check{w} \right\|_{2} < \varepsilon$.
Take any $\varepsilon>0$, $\rho>1+\varepsilon$, and $y\in \mathrm{Ball}(\check{w}, \rho-\varepsilon)$.
For $n\geq n_{\varepsilon}$, since we have
$\left\| y - {w}^{\prime}_{(n)} \right\|_{2} \leq \left\| y - \check{w} \right\|_{2} + \left\| \check{w} - {w}^{\prime}_{(n)}\right\|_{2} \leq \rho$,
we have $y \in \left\{ {w}^{\prime}_{(n)} \right\} + \mathrm{Ball}(0, \rho)$.
Thus, we have 
\[
	\mathrm{Ball}(\check{w}, 1) 
	\subset \mathrm{Ball}(\check{w}, \rho-\varepsilon) 
	\subset \left\{ {w}^{\prime}_{(n)} \right\} + \mathrm{Ball}(0, \rho)
\]
for $n\geq n_{\varepsilon}$.
Hence, 
\begin{align*}
	\bigcap_{n \geq n_\varepsilon} \left\{ \frac{1}{T}R({u}^{\prime}_{(n)}) - \mathrm{Ball}(0,\rho) \right\}
	\supset \bigcap_{n \geq n_\varepsilon} \left\{- \left\{{w}^{\prime}_{(n)}\right\} - \mathrm{Ball}(0,\rho) \right\}
	\supset - \mathrm{Ball}(\check{w}, 1).
\end{align*}
It follows from \cite[Lemma 6.6.5]{Fat99} that $(z_0, z) \neq (0,0)$, 
where we use the fact that the contingent cone of $L^2_{[0,1]}$ at ${x}^{\prime}_{(n)}$ is $L^2_{[0,1]}$, and we take a precompact sequence $\{0\}$.
Since we have $z=0$ and $z_0\in\mathbb{R}_{\geq0}$, we have $z_0 > 0$.
Then, the result follows from \eqref{eq:min_cond} and \eqref{eq:min_cond_dual}.
\end{IEEEproof}

From Theorem~\ref{thm:max-cond-convex}, we can show the discreteness property, known as the {\em bang-off-bang property}, and the optimality in Problem~\ref{prob:main_limit}, for the $L^1$ optimal control under the following assumption. 
For applications to large-scale finite-dimensional networked systems discussed in Section~\ref{sec:approximation}, we also provide a sufficient and necessary condition for this assumption to hold when the operator $A$ is compact and self-adjoint in Proposition~\ref{prop:ass_condition}. 
This characterization is useful for numerical computations (see Section~\ref{sec:example}).
These results can be seen as natural extensions of~\cite{NagQueNes16} to infinite-dimensional systems. 

\begin{assumption}\label{ass:BOB}
For $A\in\mathcal{L}(L^2_{[0,1]})$ and $b_j \in L^2_{[0,1]}$,
define a bounded linear operator $\mathcal{M}_j: \mathbb{R} \to L^2_{[0,1]}$ by $\mathcal{M}_j u_j = A b_j u_j$,
where $u_j\in\mathbb{R}$ and $j\in \left\llbracket m \right\rrbracket$.
The system $(A; \mathcal{M}_j)$ is approximately controllable\footnote{
The system $(A; \mathcal{M}_j)$ is said to be {\em approximately controllable} on $[0, T]$ if given an arbitrary $\varepsilon>0$ it is possible to steer the state from the origin to within a distance $\varepsilon$ from all points in the state space at time $T$ with a control $u\in L^2_{[0, T]}$ (see \cite[Definition 6.2.1]{CurZwa20}).} 
on $[0, T]$ for all $j\in \left\llbracket m \right\rrbracket$.
\end{assumption}

\begin{proposition}\label{prop:ass_condition}
For $A\in\mathcal{L}(L^2_{[0,1]})$ and $b \in L^2_{[0,1]}$, define a bounded linear operator $\mathcal{M}: \mathbb{R} \to L^2_{[0,1]}$ by $\mathcal{M}u = A b u$ for $u\in\mathbb{R}$.
Suppose that the operator $A$ is compact and self-adjoint, and let $\{\lambda_i\} \subset \mathbb{R}$ be eigenvalues of $A$ and $\{\phi_{i}\} \subset L^2_{[0,1]}$ be the corresponding eigenvectors that forms an orthonormal basis of $L^2_{[0,1]}$.
The system $(A; \mathcal{M})$ is approximately controllable on $[0, T]$ if and only if 
$\lambda_i \neq 0$ for all $i$, 
$\lambda_i \neq \lambda_j$ for all $i\neq j$, 
and $(b, \phi_i) \neq 0$ for all $i$.
\end{proposition}

\begin{IEEEproof}
See Appendix~\ref{sec:app_assumption}.
\end{IEEEproof}

\begin{theorem}\label{thm:discrete-convex}
Suppose that Assumption~\ref{ass:BOB} holds.
Then, any optimal solution to Problem~\ref{prob:convex_limit} takes values belonging to the set $\{0, \pm 1\}^m$ almost everywhere.
\end{theorem}

\begin{IEEEproof}
Let $\check{u}$ be any optimal solution to Problem~\ref{prob:convex_limit} and $\check{x}$ be the corresponding state.
Define
\begin{equation}\label{def:theta_jt}
	\check{\theta}_j(t) = 2\lambda\left( \check{x}(T) - x_{f}, e^{A(T-t)} b_j\right).
\end{equation}
From Theorem~\ref{thm:max-cond-convex}, for any $j\in \left\llbracket m \right\rrbracket$, we have
\begin{equation}\label{eq:opt_convex_limit}
\begin{aligned}
	\check{u}_{j}(t) 
	&\in \argmin_{v\in[-1,1]} \left\{|v_j| +  \check{\theta}_j(t) v_j  \right\} \\
	&= \begin{cases}
	\{-1\}, &\mbox{if~} \check{\theta}_j(t) > 1,\\
	\{0\}, &\mbox{if~} |\check{\theta}_j(t)| < 1,\\
	\{1\}, &\mbox{if~} \check{\theta}_j(t) < -1,\\
	[-1,0], &\mbox{if~} \check{\theta}_j(t) = 1,\\
	[0,1], &\mbox{if~} \check{\theta}_j(t) = -1,
	\end{cases}
\end{aligned}
\end{equation}
where $\check{u}_{j}(t)$ is the $j$th component of $\check{u}(t)$. 
Here, let us suppose that we have
\begin{equation}\label{eq:theta_1}
	\mu\left(\{ t\in(0,T): \check{\theta}_j(t) = 1\}\right) > 0
\end{equation}
for some $j\in \left\llbracket m \right\rrbracket$.
It follows from Lemma~\ref{lem:analytic} in Appendix~\ref{sec:app_analytic} and~\cite[Corollary 1.2.7]{Kra02} that 
for all $t \in [0,T]$, we have $\check{\theta}_j(t) = 1$ and
$\left( \check{x}(T) - x_{f}, e^{At} A b_j\right) = 0$.
Hence, 
\begin{align*}
	0 
	&= \left( \check{x}(T) - x_{f}, e^{At} \mathcal{M}_j(1) \right)\\
	&= \left( {\mathcal{M}_j}^{\prime} {e^{At}}^{\prime} \left(\check{x}(T) - x_{f}\right), 1 \right)\\
	&= {\mathcal{M}_j}^{\prime} {e^{At}}^{\prime} \left(\check{x}(T) - x_{f}\right) (1, 1),
\end{align*}
where ${\mathcal{M}_j}^{\prime}$ and ${e^{At}}^{\prime}$ denote the adjoint operator of $\mathcal{M}_j$ and $e^{At}$, respectively.
This gives 
\begin{align*}
	{\mathcal{M}_j}^{\prime} {e^{At}}^{\prime} \left(\check{x}(T) - x_{f}\right) = 0,
\end{align*}
which implies $\check{x}(T) - x_{f} = 0$ from Assumption~\ref{ass:BOB} and~\cite[Theorem 6.2.6]{CurZwa20}.
Then, we have $\check{\theta}_{j}(t) = 0$ for all $t\in[0,T]$ from the definition~\eqref{def:theta_jt}, which contradicts to \eqref{eq:theta_1}.
Hence, 
\begin{equation}\label{eq:BOB_1}
	\mu(\{ t\in(0,T): \check{\theta}_j(t) = 1\}) = 0
\end{equation}
for all $j\in \left\llbracket m \right\rrbracket$. 
Similarly, we have 
\begin{equation}\label{eq:BOB_2}
	\mu(\{ t\in(0,T): \check{\theta}_j(t) = -1\}) = 0
\end{equation}
for all $j\in \left\llbracket m \right\rrbracket$ under Assumption~\ref{ass:BOB}.
The result follows from~\eqref{eq:opt_convex_limit}, \eqref{eq:BOB_1}, and \eqref{eq:BOB_2}.
\end{IEEEproof}

Now, we are ready to show the main result of Section~\ref{subsec:analysis_1}.

\begin{theorem}\label{thm:L0-L1}
Suppose that Assumption~\ref{ass:BOB} holds.
Then, any optimal solution to Problem~\ref{prob:convex_limit} is an optimal solution to Problem~\ref{prob:main_limit}.
\end{theorem}

\begin{IEEEproof}
Let $\check{v}$ be any optimal solution to Problem~\ref{prob:convex_limit}.
Note that we have
\begin{equation}\label{eq:equivalence_1}
	\check{v}(t) \in \{0, \pm 1\}^m
\end{equation}
almost everywhere from Theorem~\ref{thm:discrete-convex}.
Note also that we have 
\begin{equation}\label{eq:equivalence_2}
	\left\|u\right\|_{1} \leq \left\|u\right\|_{0}
\end{equation}
for any $u\in\mathcal{U}$.
Hence, we have
\[
	J_0(\check{v}) = J_1(\check{v}) \leq J_1(u) \leq J_0(u)
\]
for any $u\in\mathcal{U}$,
where 
the first relation follows from \eqref{eq:equivalence_1},
the second relation follows from the optimality of $\check{v}$,
and the third relation follows from \eqref{eq:equivalence_2}.
This completes the proof.
\end{IEEEproof}

\subsection{Analysis on non-convex optimal controls}
\label{subsec:analysis_2}

We here present non-convex optimal control problems which always give a sparse optimal control.
The problems are formulated as follows, where the integrand $\psi:\mathbb{R}^m\to\mathbb{R}$ in the cost function is a given non-convex function satisfying an assumption (Assumption~\ref{ass:psi} below).

\begin{problem}[non-convex optimal control]\label{prob:nonconvex_limit}
For $A\in\mathcal{L}(L^2_{[0,1]})$, $B \in \mathscr{B}$, $x_{0}, x_{f}\in L^{2}_{[0,1]}$, $T>0$, and $\lambda>0$, 
find a control input $u$ on $[0, T]$ that solves the following:
\begin{equation*}
\begin{aligned}
  & \underset{u}{\text{minimize}}
  & & \int_{0}^{T} \psi(u(t)) dt + \lambda \left\|x(T) - x_{f}\right\|_{2}^2\\
  & \text{subject to}
  & & \dot{x}(t) = Ax(t) + Bu(t), \quad x(0) = x_{0}, \quad u \in \mathcal{U}.
\end{aligned}
\end{equation*}
\end{problem}

Here, we introduce the assumption on the function $\psi$.
The assumption is derived from the study~\cite{Ike25}, which proposed the penalty functions that give sparse optimal controls for finite-dimensional linear systems.
The assumption is satisfied by several well-known functions, such as the minimax concave penalty and the $\ell^p$ norm with $p \in (0, 1)$ (see~\cite[Remark 3]{Ike25}).

\begin{assumption}\label{ass:psi}
The function $\psi$ satisfies the following:
\begin{enumerate}
\item[(A1)] There exist continuous functions $\psi_j: \mathbb{R}\to\mathbb{R}$, $j\in\left\llbracket m \right\rrbracket$, that satisfy
	$\psi(u)=\sum_{j=1}^{m} \psi_j(u_j)$,
	where $u = \left[u_1, u_2, \dots, u_m\right]^\top$.
\item[(A2)] $\psi_j(0)=0$, 
			$|u_j| < \psi_j(u_j) \leq 1$ on $(-1, 1) \backslash \{0\}$, 
			and $\psi_j(1) = \psi_j(-1) = 1$ for all $j\in\left\llbracket m \right\rrbracket$.
\end{enumerate}
\end{assumption}

We have the following result.

\begin{theorem}\label{thm:L0-nonconvex}
Fix any function $\psi$ that satisfies Assumption~\ref{ass:psi} and consider Problem~\ref{prob:nonconvex_limit} defined by the function $\psi$.
Then, any optimal solution to Problem~\ref{prob:nonconvex_limit} is an optimal solution to Problem~\ref{prob:main_limit}.
\end{theorem}

\begin{IEEEproof}
This result follows from a similar method as in the proof of Theorem~\ref{thm:L0-L1}.
Let $\tilde{u}$ be any optimal solution to Problem~\ref{prob:nonconvex_limit} and $\tilde{x}$ be the corresponding state.
It follows from~\cite[Theorem 6.6.2]{Fat99} and \cite[Lemma 6.6.5]{Fat99} that we have
\begin{equation*}
	\tilde{u}(t) \in \argmin_{v\in[-1,1]^m} 
	\left\{\psi(v) +  2\lambda \left( \tilde{x}(T) - x_{f}, e^{A(T-t)} Bv \right) \right\}
\end{equation*}
almost everywhere.
From (A1) in Assumption~\ref{ass:psi}, this is reduced to the component-wise minimization
\begin{equation*}
	\tilde{u}_{j}(t) \in \argmin_{v_j \in [-1,1]} 
	\left\{\psi_j(v_j) + \tilde{\theta}_j(t) v_j \right\}
\end{equation*}
for all $j\in\left\llbracket m \right\rrbracket$, where 
\[
	\tilde{\theta}_j(t) 
	= 2\lambda \left( \tilde{x}(T) - x_{f}, e^{A(T-t)} b_j \right).
\]
It follows from (A2) that we have 
\begin{align*}
	\tilde{u}_{j}(t) \in 
	\begin{cases}
		\{-1\}, & \mbox{if~} \tilde{\theta}_j(t)  > 1,\\
		\{-1, 0\}, & \mbox{if~} \tilde{\theta}_j(t)  = 1,\\
		\{0\}, & \mbox{if~} |\tilde{\theta}_j(t)| < 1,\\
		\{0, 1\}, & \mbox{if~} \tilde{\theta}_j(t) = -1,\\
		\{1\}, & \mbox{if~} \tilde{\theta}_j(t) < -1
	\end{cases}
\end{align*}
almost everywhere for all $j\in \left\llbracket m \right\rrbracket$.
Hence, the control $\tilde{u}$ takes the values in $\{0, \pm1\}^m$ almost everywhere, and we have
\begin{equation}\label{eq:nonconvex_L0}
	\int_{0}^{T} \psi(\tilde{u}(t)) dt = \left\|\tilde{u}\right\|_0.
\end{equation}
Denote by $J_{\mathrm{NC}}$ the cost function of Problem~\ref{prob:nonconvex_limit}.
Then, for any $u \in \mathcal{U}$ we have
\[
	J_0(\tilde{u}) = J_{\mathrm{NC}}(\tilde{u}) \leq J_{\mathrm{NC}}(u) \leq J_0(u),
\]
where 
the first relation follows from~\eqref{eq:nonconvex_L0},
the second relation follows from the optimality of $\tilde{u}$ in Problem~\ref{prob:nonconvex_limit},
and the third relation follows from (A1) and (A2).
This completes the proof.
\end{IEEEproof}

\section{Sparse approximate control for large-scale finite-dimensional networked systems}
\label{sec:approximation}

This section investigates the sparse optimal control problem for large-scale finite-dimensional networked systems by utilizing the results in the previous section.
To represent the underlying network structures, we employ graphons that capture the overall connection pattern. 
We first review the fundamentals of graphs and graphons to set up the discussion.

\subsection{Graphs and graphons}\label{subsec:graphon}

A weighted graph consists of
a finite set of nodes $V$, 
the set of edges $E \subseteq V \times V$,
and the edge weights $a_{ij}\in\mathbb{R}$ for each edge $(i, j) \in E$.
For convenience, we set $a_{ij}=0$ if $(i, j) \not\in E$.
Throughout the paper, we consider undirected weighted graphs (i.e., $a_{ij}=a_{ji}$), and the edge weights are bounded by $a_{ij} \in [0, 1]$ for all $i, j \in V$.
The adjacency matrix of a graph with a node set $V = \left\llbracket n \right\rrbracket$ and edge weights $a_{ij}$ is defined by a matrix $[a_{ij}]\in\mathbb{R}^{n\times n}$.

Let $\mathcal{W}$ denote the space of all bounded measurable functions $W: [0,1]^2 \to \mathbb{R}$ such that
$W(\alpha, \beta) = W(\beta, \alpha)$ for all $\alpha, \beta \in [0,1]$.
Define sets
\begin{align*}
	&\mathcal{W}_0 = \left\{W\in\mathcal{W}: 
		0 \leq W(\alpha, \beta) \leq 1, \forall (\alpha, \beta) \in [0,1]^2\right\},\\
	&\mathcal{W}_1 = \left\{W\in\mathcal{W}: 
		|W(\alpha, \beta)| \leq 1, \forall (\alpha, \beta) \in [0,1]^2\right\}.
\end{align*}
The elements of $\mathcal{W}_0$ are called {\em graphons}.
For any graphon $W$, the value $W(\alpha, \beta)$ can be interpreted as the edge weight between nodes $\alpha$ and $\beta$, and the interval $[0,1]$ represents the node set.
Graphons generalize weighted graphs in the following sense:
For a graph $G$ with a node set $V = \left\llbracket n \right\rrbracket$ and edge weights $a_{ij}$,
we can obtain a graphon $W$ by setting $W(\alpha, \beta) = a_{ij}$ if $(\alpha, \beta) \in P_i \times P_j$,
where 
$\{P_{i}\}$ denotes the set of uniformly partitioned subintervals of $[0,1]$, i.e.,
$P_i = \left[\frac{i-1}{n}, \frac{i}{n}\right)$ for $i \in \llbracket n-1 \rrbracket$ and $P_n = \left[\frac{n-1}{n}, 1\right]$.
For a graph $G$, the graphon obtained from $G$ via the above construction is denoted by $W_{G}$.
(In general, a graphon that is constant over each $P_i \times P_j$ is called {\em a step graphon}.)

Every graphon $W\in\mathcal{W}_0$ defines a bounded linear operator $T_W \in \mathcal{L}(L^{2}_{[0,1]})$ by
\begin{equation}\label{eq:operator_graphon}
	(T_W x)(\alpha) = \int_{0}^{1} W(\alpha, \beta) x(\beta) d\beta.
\end{equation}
(Indeed, we have $\left\|T_{W}\right\|_{\mathrm{op}}\leq 1$, since $\left\|T_{W}x\right\|_{2} \leq \left\|x\right\|_2$ by Cauchy-Schwarz inequality.)
Intuitively, the operator $T_W$ plays the same role as the adjacency matrix in finite graphs.
Let $\mathscr{A}$ denote the set of all the operators, i.e., 
$\mathscr{A} = \left\{ T_W: W \in \mathcal{W}_0 \right\}$.
The {\em cut norm} on the linear space $\mathcal{W}$ is defined by
\[
	\left\|W\right\|_{\square} = \sup_{S, T \subset [0,1]} \left|\int_{S \times T} W(\alpha, \beta) d\alpha d\beta\right|,
\]
where the supremum is taken over all measurable subsets $S$ and $T$.
For every $W\in\mathcal{W}_1$, we have the following relationship (see~\cite[Lemma E.6]{Jan13} and \cite[Lemma 8.11]{lov12}):
\begin{align}\label{ineq:norm_cut_op}
	\left\|W\right\|_{\square} \leq \left\|T_W\right\|_{\mathrm{op}} \leq 2\sqrt{2} \left\|W\right\|_{\square}^{1/2},
\end{align}
where $T_W$ is defined by~\eqref{eq:operator_graphon}.
From this,
for a sequence of graphons $\{W^{[n]}\}$ and a graphon $W\in\mathcal{W}_0$,
it holds that $\left\|T_{W^{[n]}} - T_W\right\|_{\mathrm{op}} \to 0$ if and only if $\left\|W^{[n]}-W\right\|_{\square} \to 0$.

To evaluate the similarity between two graphs, the {\em homomorphism density} has been used.
This quantity intuitively indicates how likely a small graph appears as a pattern inside a larger graph
(see~\cite[Chapter~5]{lov12} for the precise definition).
This notion has been employed as a criterion for the convergence of graph sequences. 
More precisely, a sequence of graphs $\{G^{[n]}\}$ is said to be {\em convergent} if the homomorphism density $t(F,G^{[n]})$ of $F$ in $G^{[n]}$ converges for any finite simple graph $F$.
Homomorphism densities in graphs can be naturally extended to homomorphism densities in graphons~\cite[Chapter~7]{lov12}, and it is known that
for any convergent sequence of graphs $\{G^{[n]}\}$, there exists a graphon $W\in\mathcal{W}_0$ such that $t(F,G^{[n]}) \to t(F, W)$ for all finite simple graph $F$~\cite[Theorem 3.8]{borChaLov12}.
This graphon $W$ is called the {\em limit} of the graph sequence $\{G^{[n]}\}$, and this is denoted by $G^{[n]} \to W$.
Furthermore, it is also known that if $G^{[n]} \to W$ for a graphon $W \in\mathcal{W}_0$, then the graphs can be labeled so that $\left\|W_{G^{[n]}} - W\right\|_{\square}\to0$~\cite[Lemma 5.3]{borChaLov12}.
Conversely, if $\left\|W_{G^{[n]}} - W\right\|_{\square}\to0$ for a graphon $W\in\mathcal{W}_0$, then we have $G^{[n]} \to W$.
Then, this section considers finite-dimensional systems whose network structures are defined by a graph $G^{[n]}$ such that $W_{G^{[n]}}$ converges to a graphon $W$ with respect to the cut norm (see Assumption~\ref{ass:convergence}). 

\subsection{Problem formulation}\label{subsec:formulation_2}

We here formulate our sparse optimal control problem for finite-dimensional systems.
We consider a finite-dimensional system consisting of $n$ nodes whose dynamics are defined by
\[
	\dot{\mathrm{x}}_{i}^{[n]}(t) = \frac{1}{n} \sum_{j=1}^{n} \mathcal{A}_{ij}^{[n]} \mathrm{x}_{i}^{[n]}(t) + \sum_{j=1}^{m} \mathcal{B}_{j,i}^{[n]}u_{j}(t),
\]
where
$\mathcal{A}_{ij}^{[n]}\in[0,1]$ is an edge weight of the underlying undirected graph $G^{[n]}$ between nodes,
$\mathcal{B}_{j,i}^{[n]}\in \mathbb{R}$ is the $i$th element of $\mathcal{B}_j^{[n]}\in \mathbb{R}^{n}$,
and $\mathrm{x}_{i}^{[n]} (t) \in \mathbb{R}$ is the state of node $i$ at time $t$.
Define
$\mathcal{A}^{[n]} = \left[\mathcal{A}_{ij}^{[n]}\right]\in\mathbb{R}^{n\times n}$,
$\mathcal{B}^{[n]} = \left[\mathcal{B}_1^{[n]}, \mathcal{B}_2^{[n]}, \dots, \mathcal{B}_m^{[n]}\right] \in\mathbb{R}^{n\times m}$,
and $\mathrm{x}^{[n]}(t) = \left[\mathrm{x}_{1}^{[n]}(t), \mathrm{x}_{2}^{[n]}(t), \dots, \mathrm{x}_{n}^{[n]}(t)\right]^\top \in \mathbb{R}^{n}$ with an initial state $\mathrm{x}_{0}^{[n]}\in\mathbb{R}^{n}$.
Then, the overall system is described by
\begin{equation}\label{eq:system_finite}
	\dot{\mathrm{x}}^{[n]}(t) = \mathcal{A}^{[n]} \circ \mathrm{x}^{[n]}(t) + \mathcal{B}^{[n]}u(t)
\end{equation}
with $\mathrm{x}^{[n]}(0) = \mathrm{x}_{0}^{[n]}$,
where $\circ$ denotes the averaging operator defined by $\mathcal{A}^{[n]} \circ \mathrm{x} = \frac{1}{n} \mathcal{A}^{[n]} \mathrm{x} $ for $\mathrm{x} \in \mathbb{R}^{n}$.

Here, define 
operators $A^{[n]}: L^2_{[0,1]} \to L^2_{[0,1]}$, $B^{[n]}: \mathbb{R}^m \to L^2_{[0,1]}$, and $M^{[n]}: \mathbb{R}^n\to L^2_{[0, 1]}$ by
\begin{align}
	&A^{[n]} = T_{W_{G^{[n]}}}, \label{eq:op_An}\\
	&\left[B^{[n]} u\right](\alpha) = \sum_{j=1}^{m} b_{j}^{[n]}(\alpha) u_j, \quad \alpha \in [0, 1], \label{eq:op_Bn}\\
	&\left[M^{[n]}(\mathrm{x})\right](\alpha) = \sum_{i=1}^{n} \mathbf{1}_{P_i}(\alpha) \mathrm{x}_{i}, \quad \alpha \in [0, 1], \label{eq:map_M}
\end{align}
where 
$T_{W_{G^{[n]}}}$, $\{P_{i}\}$, and $W_{G^{[n]}}$ are defined in Section~\ref{subsec:graphon} 
(i.e., the operator \eqref{eq:operator_graphon} for the graphon $W_{G^{[n]}}$,
subintervals obtained by uniformly partitioning the interval $[0,1]$, 
and the step graphon constructed from $G^{[n]}$), 
and $b_{j}^{[n]}$ is a piecewise constant function defined by
\begin{align*}
	&b_{j}^{[n]} (\alpha) = \left[M^{[n]} (\mathcal{B}_{j}^{[n]})\right](\alpha), \quad \alpha \in [0, 1], \quad j\in\left\llbracket m \right\rrbracket.
\end{align*}
Then, we have the following lemma.
Note that in this paper state variables are denoted by the upright (roman) symbol (e.g. ${\rm x}^{[n]}$ and ${\rm x}_0^{[n]}$) for finite-dimensional systems, while they are denoted by the italic symbol (e.g. $x^{[n]}$ and $x_0^{[n]}$; see below) for infinite-dimensional systems to avoid confusion.

\begin{lemma}\label{lem:finite_infinite_system}
Consider an infinite-dimensional system
\begin{equation}\label{eq:system_finite_ver2}
	\dot{x}^{[n]}(t) = A^{[n]} x^{[n]}(t) + B^{[n]}u(t)
\end{equation}
with the initial state $x_{0}^{[n]} = M^{[n]}(\mathrm{x}_{0}^{[n]})$,
where 
$A^{[n]}$, $B^{[n]}$, $M^{[n]}$ are defined by \eqref{eq:op_An}, \eqref{eq:op_Bn}, and \eqref{eq:map_M}, respectively.
Then, the state trajectories of the system~\eqref{eq:system_finite} correspond one-to-one to that of the system \eqref{eq:system_finite_ver2} under the mapping $M^{[n]}$.
\end{lemma}
\begin{IEEEproof}
See \cite[Lemma 3]{GaoCai20}.
\end{IEEEproof}

Based on Lemma~\ref{lem:finite_infinite_system}, the representation~\eqref{eq:system_finite_ver2} for finite-dimensional systems is considered in the subsequent discussion.
We denote by $(A^{[n]}; B^{[n]})$ or $(A^{[n]}; B^{[n]}; x_0^{[n]})$ the system~\eqref{eq:system_finite_ver2}.

Let the target state of finite-dimensional systems be denoted by $\mathrm{x}_f^{[n]} \in \mathbb{R}^n$. 
For the system \eqref{eq:system_finite_ver2}, the target state is expressed as $x_f^{[n]} = M^{[n]}(\mathrm{x}_f^{[n]}) \in L^2_{[0,1]}$.  
Accordingly, the cost function is defined by
\begin{equation*}
	J_0^{[n]}(u) = \left\|u\right\|_{0} + \lambda \left\| x^{[n]}(T) - x_f^{[n]}\right\|_2^2.
\end{equation*}
Now, we are ready to formulate the sparse optimal control problem for the finite-dimensional system.

\begin{problem}\label{prob:sparse_finite}
	Find a control $u \in \mathcal{U}$ that minimizes $J_0^{[n]}(u)$ for the system $(A^{[n]}; B^{[n]}; x_0^{[n]})$.
\end{problem}

In large-scale systems, it is difficult to obtain the exact network structure, and control designs that depend on the exact system size may lead to undesirable consumption of computational resources.
To address these issues, we consider an approach that computes approximate solutions based on the underlying connection pattern of the system, instead of directly solving the problem (for which the optimal feedback law is studied in~\cite{IkeKas_ECC19}).
In other words, we consider the scenario in which $W_{G^{[n]}}$, $B^{[n]}$, $x_{0}^{[n]}$, and $x_f^{[n]}$ in the system~\eqref{eq:system_finite_ver2} converge to known limit objects.
More precisely, we put the following assumption:

\begin{assumption}\label{ass:convergence}
For given $W\in\mathcal{W}_0$, $b_1, b_2, \dots, b_m \in L^2_{[0, 1]}$, and $x_{0}, x_{f} \in L^2_{[0, 1]}$, the following conditions hold:
\begin{enumerate}
    \item $ \displaystyle \lim_{n \to \infty} \left\| W - W_{G^{[n]}} \right\|_{\square} = 0$.
    \item $ \displaystyle \lim_{n \to \infty} \left\| b_{j} - b_{j}^{[n]} \right\|_{2} = 0$ for all $j\in\left\llbracket m \right\rrbracket$.
    \item $ \displaystyle \lim_{n \to \infty} \left\| x_{0} - x_{0}^{[n]} \right\|_{2} = \displaystyle \lim_{n \to \infty} \left\| x_{f} - x_{f}^{[n]} \right\|_{2} = 0$.
\end{enumerate}
\end{assumption}

We call the system~\eqref{eq:system} defined by the limits $W$, $\{b_j\}$, and $x_0$, with $A = T_W$, $Bu = \sum_{j=1}^{m}b_j u_j$, and $x(0) = x_0$, {\em the limit graphon system}. 
Let $\bar{u}$ be an optimal solution to Problem~\ref{prob:main_limit} for the given system~\eqref{eq:system} and $\bar{x}$ be the corresponding state.
Also, let $\bar{u}^{[n]}$ be an optimal solution to Problem~\ref{prob:sparse_finite} for the system $(A^{[n]}; B^{[n]}; x_{0}^{[n]})$ and $\bar{x}^{[n]}$ be the corresponding state.
Note that for the computation of $\bar{u}$, Theorem~\ref{thm:L0-L1} and Theorem~\ref{thm:L0-nonconvex} in the previous section are useful. 
Also, since the operator $A=T_W$ is self-adjoint~\cite[Example A.3.62]{CurZwa20} and compact~\cite[Theorem A.3.25]{CurZwa20}, Proposition~\ref{prop:ass_condition} is applicable to check the sparsity of the $L^1$ optimal control.
In Section~\ref{subsec:analysis_approximation}, we show that the optimal control $\bar{u}^{[n]}$ can be approximated by $\bar{u}$.


\subsection{Analysis}\label{subsec:analysis_approximation}

We here show that the performance of the sparse optimal control $\bar{u}$ for the limit graphon system $(A; B; x_0)$ becomes arbitrarily close to that of the optimal control $\bar{u}^{[n]}$ to Problem~\ref{prob:sparse_finite} as the number of nodes $n$ tends to infinity.
To this end, we prepare a lemma.

\begin{lemma}\label{lem:app_limit}
Suppose that Assumption~\ref{ass:convergence} holds. 
Then, we have
\begin{align}
	&\lim_{n\to\infty} J_0^{[n]}(\bar{u}) = J_0(\bar{u}),\label{eq:approx_lemma_1}\\
	&\lim_{n\to\infty} \left| J_0^{[n]}(\bar{u}^{[n]}) - J_0(\bar{u}^{[n]}) \right| = 0.\label{eq:approx_lemma_2}
\end{align}
\end{lemma}
\begin{IEEEproof}
For \eqref{eq:approx_lemma_1}, it is enough to show 
\begin{align}
	\lim_{n\to\infty} \left\| {x}^{[n]}(x_0^{[n]}; T; \bar{u}) - x_f^{[n]} \right\|_{2}^2 
	= \left\| \bar{x}(T) - x_{f} \right\|_{2}^2, \label{eq:statement1}
\end{align}
where
${x}^{[n]}(x_0^{[n]}; T; \bar{u}) $ is the state of the system $(A^{[n]}; B^{[n]}; x_{0}^{[n]})$ at time $T$ when the control $\bar{u}$ is applied.
Note that we have
\begin{align*}
	\left\|{x}^{[n]}(x_0^{[n]}; T; \bar{u}) - x_{f}^{[n]}\right\|_{2}^2 - \left\| \bar{x}(T) - x_{f}\right\|_{2}^2 
	& = \left( \left\|\pi_1^{[n]}\right\|_{2}^{2} - \left\|\pi_3\right\|_{2}^{2} \right)
		+ 2 \left( \left(\pi_1^{[n]}, \pi_2^{[n]}\right) - \left(\pi_3, \pi_4\right) \right)\\
	&\mspace{20mu} 
		+ \left( \left\|\pi_2^{[n]}\right\|_{2}^{2} - \left\|\pi_4\right\|_{2}^{2} \right)
		+ 2 \left(\left(\pi_3, x_f\right) - \left(\pi_1^{[n]}, x_f^{[n]} \right)\right)\\
	&\mspace{20mu} 
		+ 2 \left(\left(\pi_4, x_f \right) - \left(\pi_2^{[n]}, x_f^{[n]} \right)\right)
		+ \left( \left\|x_f^{[n]}\right\|_{2}^{2} - \left\|x_f\right\|_{2}^{2} \right)
\end{align*}
where
\begin{equation}\label{def:pi}
\begin{aligned}
	&\pi_1^{[n]} = e^{A^{[n]}T}x_{0}^{[n]},\quad 
		\pi_2^{[n]} = \int_{0}^{T}e^{A^{[n]}(T-t)}B^{[n]}\bar{u}(t)dt,\\
	&\pi_3 = e^{AT}x_{0},\quad
		\pi_4 = \int_{0}^{T}e^{A(T-t)}B\bar{u}(t)dt.	
\end{aligned}
\end{equation}
From Lemma~\ref{lem:convergence} in Appendix~\ref{sec:app_convergence}, we obtain \eqref{eq:statement1}.

We next show \eqref{eq:approx_lemma_2}.
Define
\begin{equation}\label{def:pi_statement3}
\begin{aligned}
	&\pi_5^{[n]} = \int_{0}^{T}e^{A^{[n]}(T-t)}B^{[n]}\bar{u}^{[n]}(t) dt, \\
	&\pi_6^{[n]} = \int_{0}^{T}e^{A(T-t)}B\bar{u}^{[n]}(t) dt.	
\end{aligned}
\end{equation}
For any $n\in\mathbb{N}$, we have
\begin{align*}
	\frac{1}{\lambda}\left( J_0^{[n]}(\bar{u}^{[n]}) - J_0(\bar{u}^{[n]}) \right)
	&= \left\|\bar{x}^{[n]}(T) - x_{f}^{[n]} \right\|_{2}^{2} 
		- \left\|x(x_{0}; T; \bar{u}^{[n]}) - x_{f} \right\|_{2}^{2} \\
	& = \left( \left\|\pi_1^{[n]}\right\|_{2}^{2} - \left\|\pi_3\right\|_{2}^{2} \right)
		+ 2 \left( \left(\pi_1^{[n]}, \pi_5^{[n]}\right) - \left(\pi_3, \pi_6^{[n]}\right) \right)\\
	&\mspace{20mu} 
		+ \left( \left\|\pi_5^{[n]}\right\|_{2}^{2} - \left\|\pi_6^{[n]}\right\|_{2}^{2} \right)
		+ 2 \left( \left(\pi_3, x_f \right) - \left(\pi_1^{[n]}, x_f^{[n]} \right)\right)\\
	&\mspace{20mu} 
		+ 2 \left( \left(\pi_6^{[n]}, x_f \right) - \left(\pi_5^{[n]}, x_f^{[n]} \right)\right)
		+ \left( \left\|x_f^{[n]}\right\|_{2}^{2} - \left\|x_f\right\|_{2}^{2} \right)
\end{align*}
where 
${x}(x_0; T; \bar{u}^{[n]})$ is the state of the system $(A; B; x_{0})$ at time $T$ when the control $\bar{u}^{[n]}$ is applied.
From Lemma~\ref{lem:convergence} in Appendix~\ref{sec:app_convergence}, we obtain \eqref{eq:approx_lemma_2}, and the proof is completed.
\end{IEEEproof}

Now, we present the main result of this section. 
The equation \eqref{result:convergence_1} in the result shows that the control $\bar{u}$ approximates the optimal solution $\bar{u}^{[n]}$ for the finite-dimensional system for large $n$. 
Furthermore, \eqref{result:convergence_2} shows that the optimal value for the finite-dimensional system converges to that for the limit graphon system, 
while \eqref{result:convergence_3} shows that the control $\bar{u}^{[n]}$ approximates the optimal solution $\bar{u}$ for the limit graphon system.

\begin{theorem}\label{thm:app_convergence}
Suppose that Assumption~\ref{ass:convergence} holds.
Then, we have
\begin{align}
	&\lim_{n\to\infty} \left|J_0^{[n]}(\bar{u}) - J_0^{[n]}(\bar{u}^{[n]}) \right| = 0, \label{result:convergence_1}\\ 
	&\lim_{n\to\infty} J_0^{[n]}(\bar{u}^{[n]}) = J_0(\bar{u}), \label{result:convergence_2}\\ 
	&\lim_{n\to\infty} J_0(\bar{u}^{[n]}) = J_0(\bar{u}). \label{result:convergence_3}
\end{align}
\end{theorem}
\begin{IEEEproof}
Take any $\varepsilon>0$.
From Lemma~\ref{lem:app_limit}, there exists a number $N_{\varepsilon}\in\mathbb{N}$ such that we have
\begin{equation}\label{eq:app_limit}
	\left| J_0^{[n]}(\bar{u}) - J_0(\bar{u}) \right| < \frac{\varepsilon}{2}
\end{equation}
and
\begin{equation}\label{eq:cost_finite}
	\left| J_0^{[n]}(\bar{u}^{[n]}) - J_0(\bar{u}^{[n]}) \right| < \frac{\varepsilon}{2}
\end{equation}
for all $n \geq N_{\varepsilon}$.
Hence, for any $n \geq N_{\varepsilon}$, we have
\begin{align*}
	J_0^{[n]}(\bar{u}^{[n]}) 
	\leq J_0^{[n]}(\bar{u})
		< J_0(\bar{u}) + \frac{\varepsilon}{2}
	\leq  J_0(\bar{u}^{[n]}) + \frac{\varepsilon}{2}
	<  J_0^{[n]}(\bar{u}^{[n]}) + \varepsilon,
\end{align*}
where 
the first relation follows from the optimality of $\bar{u}^{[n]}$,
the second relation follows from \eqref{eq:app_limit},
the third relation follows from the optimality of $\bar{u}$,
and the fourth relation follows from \eqref{eq:cost_finite}.
This gives
\begin{align}
	&0 \leq J_0^{[n]}(\bar{u}) - J_0^{[n]}(\bar{u}^{[n]})  < \varepsilon, \label{eq:convergence_1}\\
	& \left| J_0(\bar{u}) - J_0^{[n]}(\bar{u}^{[n]}) \right| < \frac{\varepsilon}{2},\label{eq:convergence_2}
\end{align}
and 
\begin{align}
\begin{split}
\label{eq:convergence_3}
	0 
	&\leq J_0(\bar{u}^{[n]}) - J_0(\bar{u})\\
	&= \left( J_0(\bar{u}^{[n]}) - J_0^{[n]}(\bar{u}^{[n]}) \right) - \left( J_0(\bar{u}) - J_0^{[n]}(\bar{u}^{[n]}) \right)\\
	&< \frac{\varepsilon}{2} + \frac{\varepsilon}{2} = \varepsilon.
\end{split}
\end{align}
The results  \eqref{result:convergence_1}, \eqref{result:convergence_2}, and \eqref{result:convergence_3} follow from \eqref{eq:convergence_1}, \eqref{eq:convergence_2}, and \eqref{eq:convergence_3}, respectively.
\end{IEEEproof}

\section{Example}\label{sec:example}

This section illustrates the proposed method through three examples. 
Sections~\ref{subsec:example_1},~\ref{subsec:example_2}, and~\ref{subsec:example_3} focus on 
the $L^1$ optimization, 
the non-convex optimization,
and the approximation for finite-dimensional systems,
respectively.

\subsection{Example 1}\label{subsec:example_1}

\begin{figure}[t]
	\centering
	\includegraphics[width=0.8\linewidth]{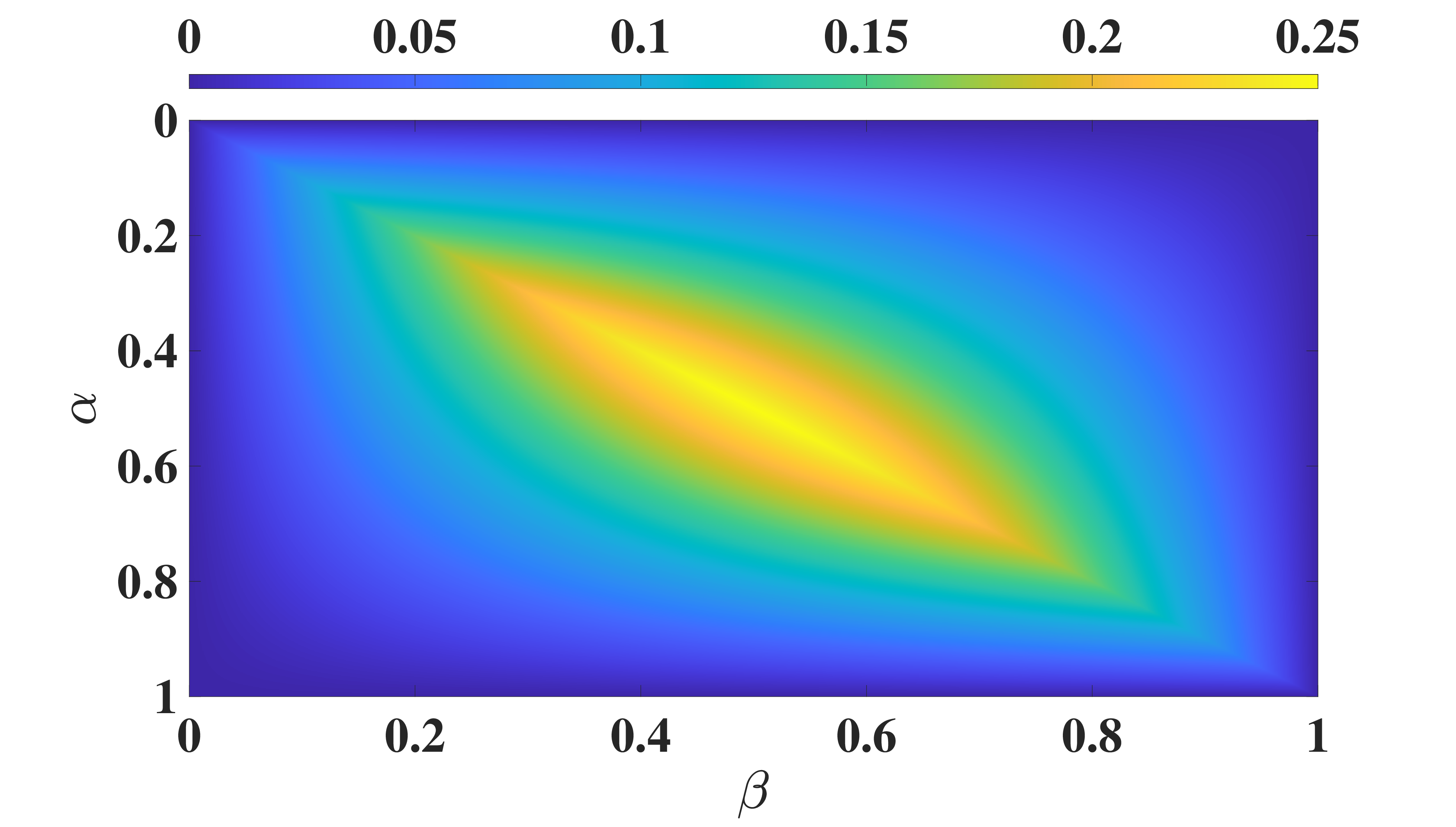}
	\caption{The graphon $W$ in Example~1.}
	\label{fig:graphon}
\end{figure}

We consider the system~\eqref{eq:system} with the operator $A=T_W$ defined by a graphon 
\begin{align*}
	W(\alpha, \beta) = 
	\begin{cases}
		(1-\beta) \alpha, &\mbox{for~} 0 \leq \alpha \leq \beta \leq 1,\\
		(1-\alpha) \beta, &\mbox{for~} 0 \leq \beta \leq \alpha \leq 1,
	\end{cases}
\end{align*}
and the operator $B$ defined by
\begin{align*}
	b_1(\alpha) = {\bf 1}_{\Omega_1}(\alpha),\quad
		b_2(\alpha) = {\bf 1}_{\Omega_2}(\alpha),\quad
	\Omega_1 = [0,0.5], \quad
		\Omega_2 = (0.5, 0.8].
\end{align*}
The graphon $W$ is shown in Figure~\ref{fig:graphon}.
The operator $A$ has eigenvalues $\frac{1}{k^2 \pi^2}$ and eigenvectors $\sqrt{2}\sin(k\pi\alpha)$, where $k=1,2,\dots$ (see~\cite[Example A.4.20]{CurZwa20}).
Hence, this example satisfies Assumption~\ref{ass:BOB} from Proposition~\ref{prop:ass_condition},
and a sparse optimal control can be obtained by solving Problem~\ref{prob:convex_limit}, as established in Theorem~\ref{thm:L0-L1}.
For the numerical optimization, we used the package CVX with MATLAB~\cite{cvx}.

\begin{figure}[t]
	\centering
	\includegraphics[width=0.8\linewidth]{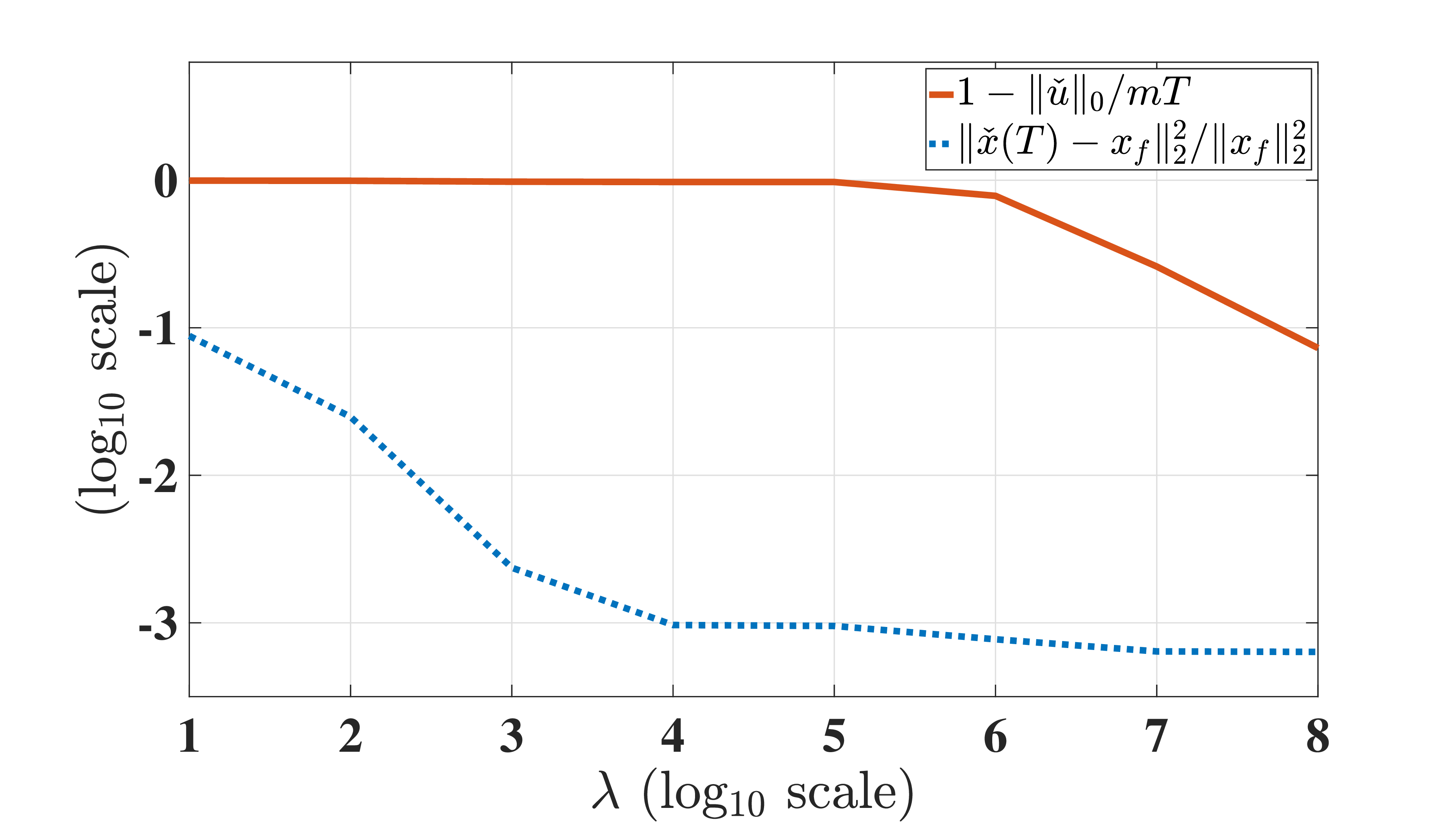}
	\caption{The effect of the parameter $\lambda$ on the sparsity rate $1 - \frac{\left\|\check{u}\right\|_{0}}{mT}$ (solid line) and the terminal error rate $\frac{\left\| \check{x}(T) - x_f\right\|_{2}^{2}}{\left\|x_f\right\|_2^2}$
 (dotted line) in Example~1.}
	\label{fig:trade_off}
\end{figure}

\begin{figure}[t]
	\centering
	\includegraphics[width=0.49\linewidth]{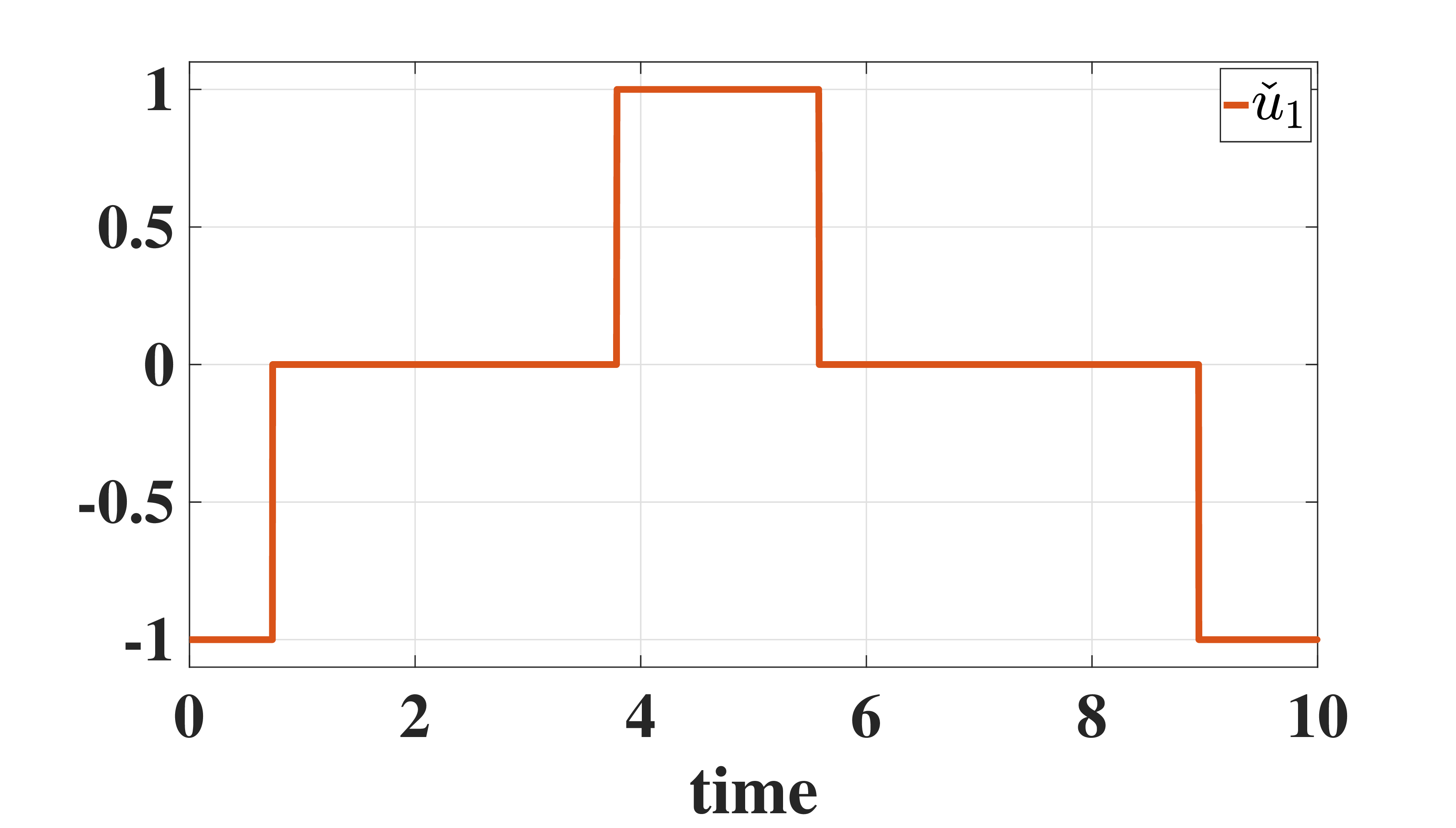}
	\includegraphics[width=0.49\linewidth]{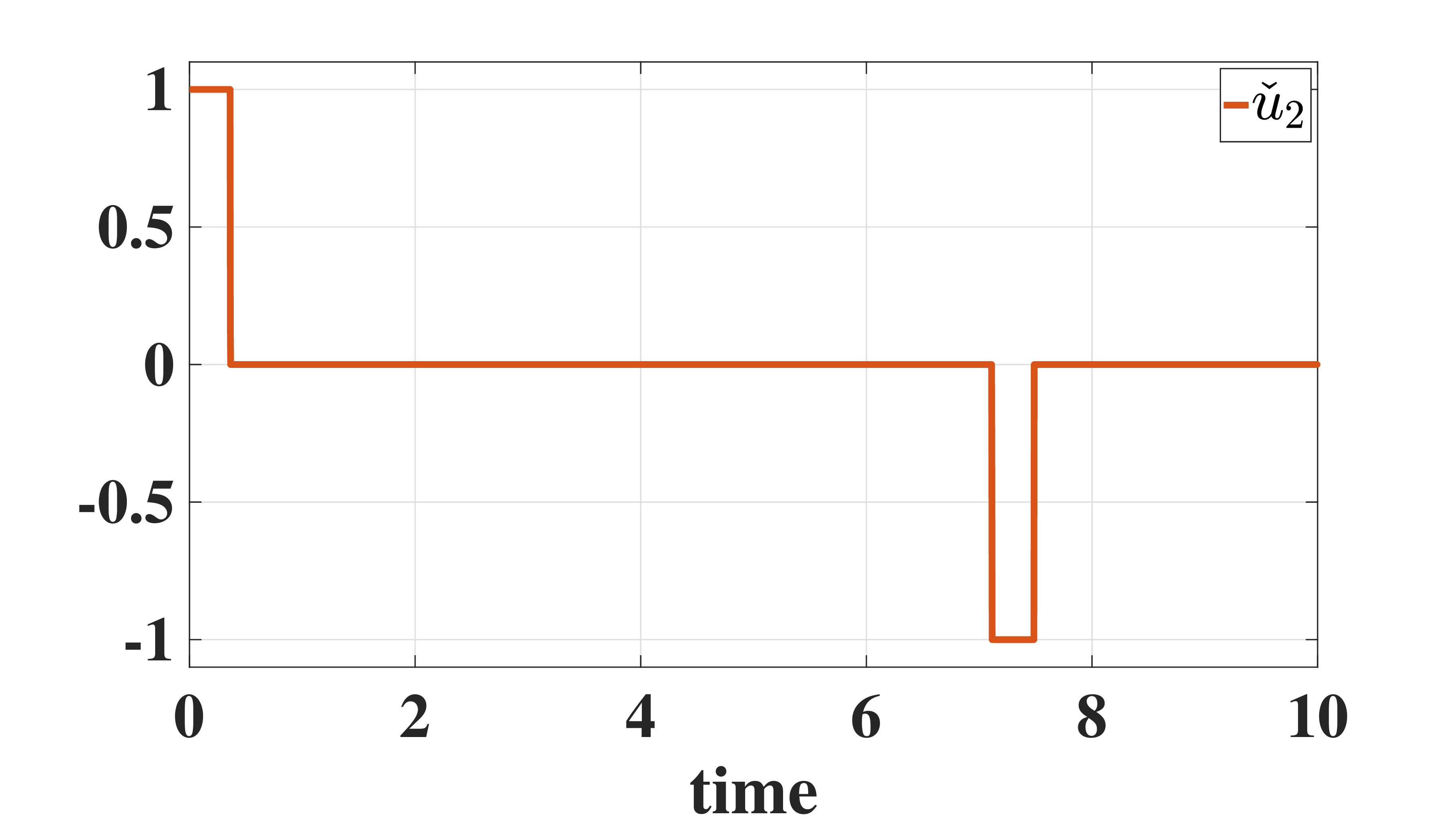}
	\caption{The $L^1$ optimal control $\check{u}$ in Example~1.}
	\label{fig:control}
\end{figure}

\begin{figure}[t]
	\centering
	\includegraphics[width=0.8\linewidth]{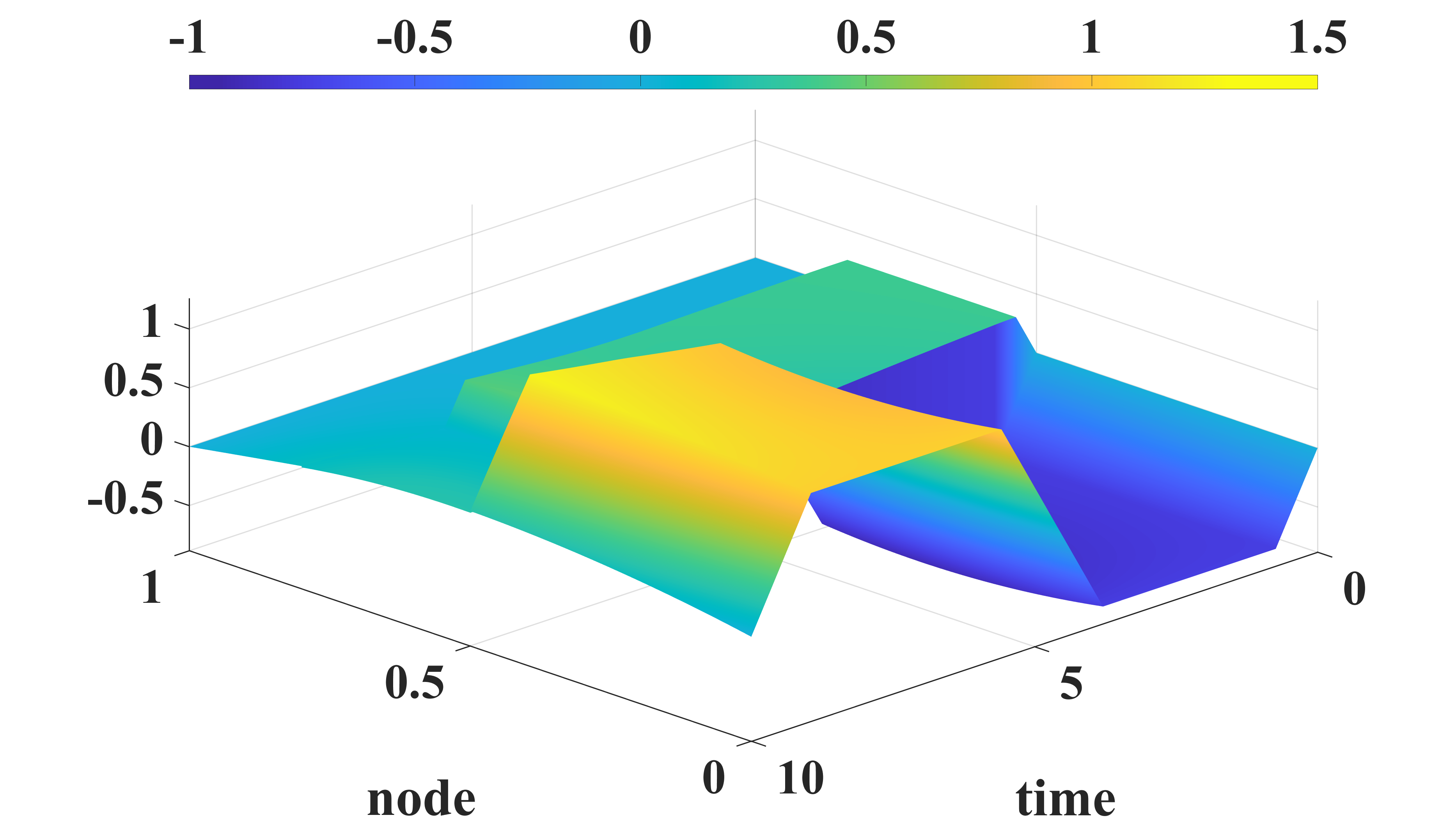}
	\caption{The state trajectories $\check{x}$ in Example~1.}
	\label{fig:state}
\end{figure}

For the sparse optimal control problem with the parameters 
$x_0 = 0$, $x_f(\alpha) = - \alpha(\alpha  - 1)$, and $T = 10$,
Figure~\ref{fig:trade_off} illustrates the effect of the regularization parameter $\lambda$ on the sparsity rate $1-\frac{\left\|\check{u}\right\|_{0}}{mT}$ and the terminal error rate $\frac{\left\| \check{x}(T) - x_f\right\|_{2}^{2}}{\left\|x_f\right\|_2^2}$, where $\check{u}$ is the obtained $L^1$ optimal control and $\check{x}$ is the corresponding state.
As $\lambda$ increases, the sparsity of the optimal control decreases, while the resulting state approaches the target state more closely.
In this example, when $\lambda$ exceeds $10^6$, the improvement in the terminal error becomes marginal, while the sparsity of the control deteriorates significantly.  
Figures~\ref{fig:control} and~\ref{fig:state} show the resulting $L^1$ optimal control and corresponding state trajectories for $\lambda = 10^6$.
It can be confirmed that the control input is certainly sparse, with a sparsity rate $1 - \frac{\left\| \check{u} \right\|_{0}}{mT} = 0.7837$, and its values are confined to $\{0, \pm1\}^2$ as established in Theorem~\ref{thm:discrete-convex}.

\subsection{Example 2}\label{subsec:example_2}

We consider the system~\eqref{eq:system} with the operator $A=T_W$ defined by a graphon $W(\alpha, \beta) = 1$
and the operator $B$ defined by
\begin{align*}
	&b_1(\alpha) = 3\alpha {\bf 1}_{\Omega_1}(\alpha) + 3(\alpha -1) {\bf 1}_{\Omega_2}(\alpha), \\
	&b_2(\alpha) = -6(\alpha - 0.5) {\bf 1}_{\Omega_3}(\alpha), \\
	&\Omega_1 = \left[0, \frac{1}{3} \right], 
		\quad \Omega_2 = \left[\frac{2}{3}, 1\right],
		\quad \Omega_3 = \left(\frac{1}{3}, \frac{2}{3}\right).
\end{align*}
For the graphon $W$, the associated operator $A = T_W$ has a zero eigenvalue (indeed, $Ab_1 = 0$).
Therefore, Assumption~\ref{ass:BOB} is not satisfied, and it is not guaranteed that the $L^1$ optimization gives a sparse optimal control.
In such cases, we can instead employ the non-convex optimal control approach, as established in Theorem~\ref{thm:L0-nonconvex}.

\begin{figure}[t]
	\centering
	\includegraphics[width=0.49\linewidth]{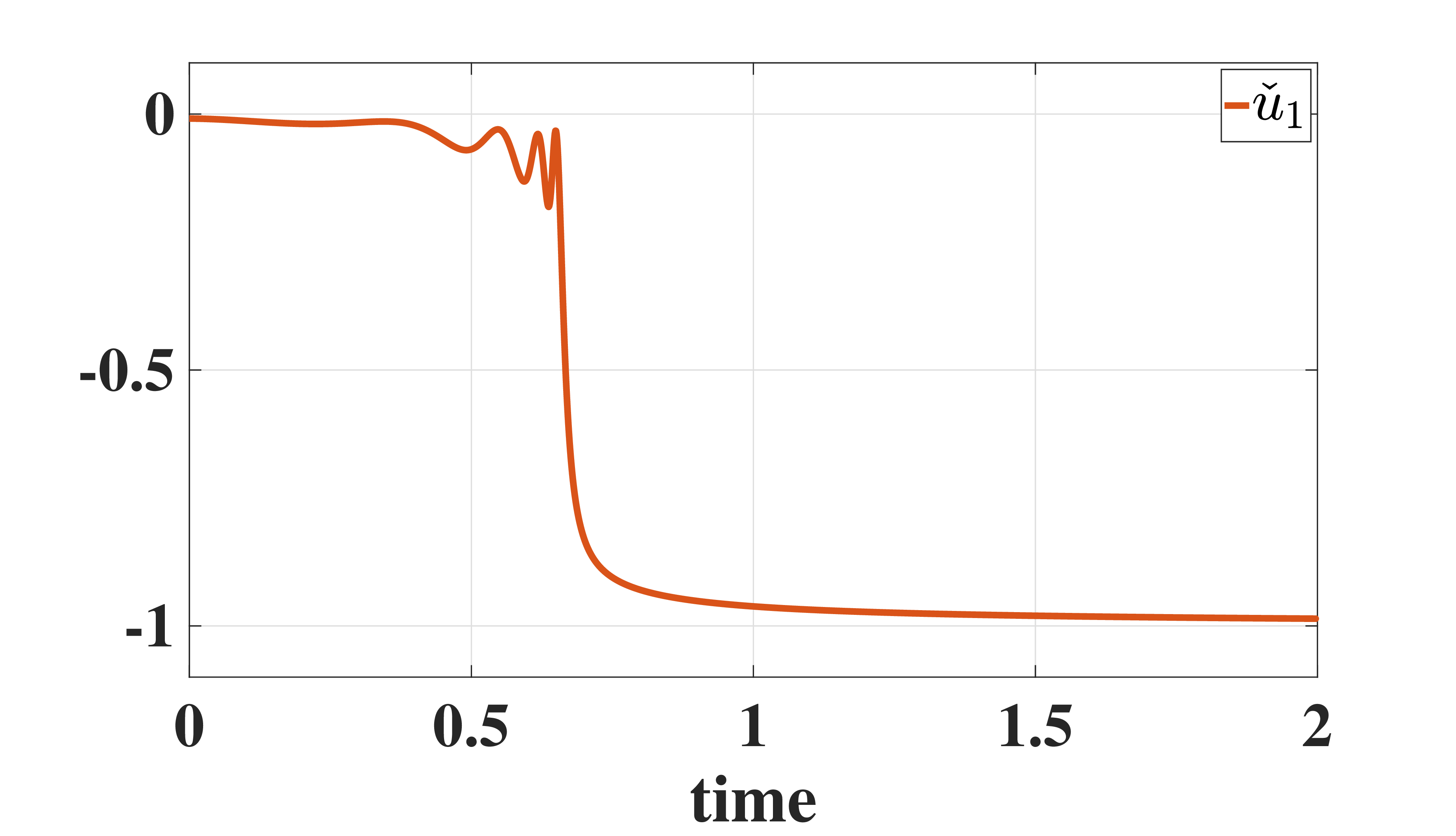}
	\includegraphics[width=0.49\linewidth]{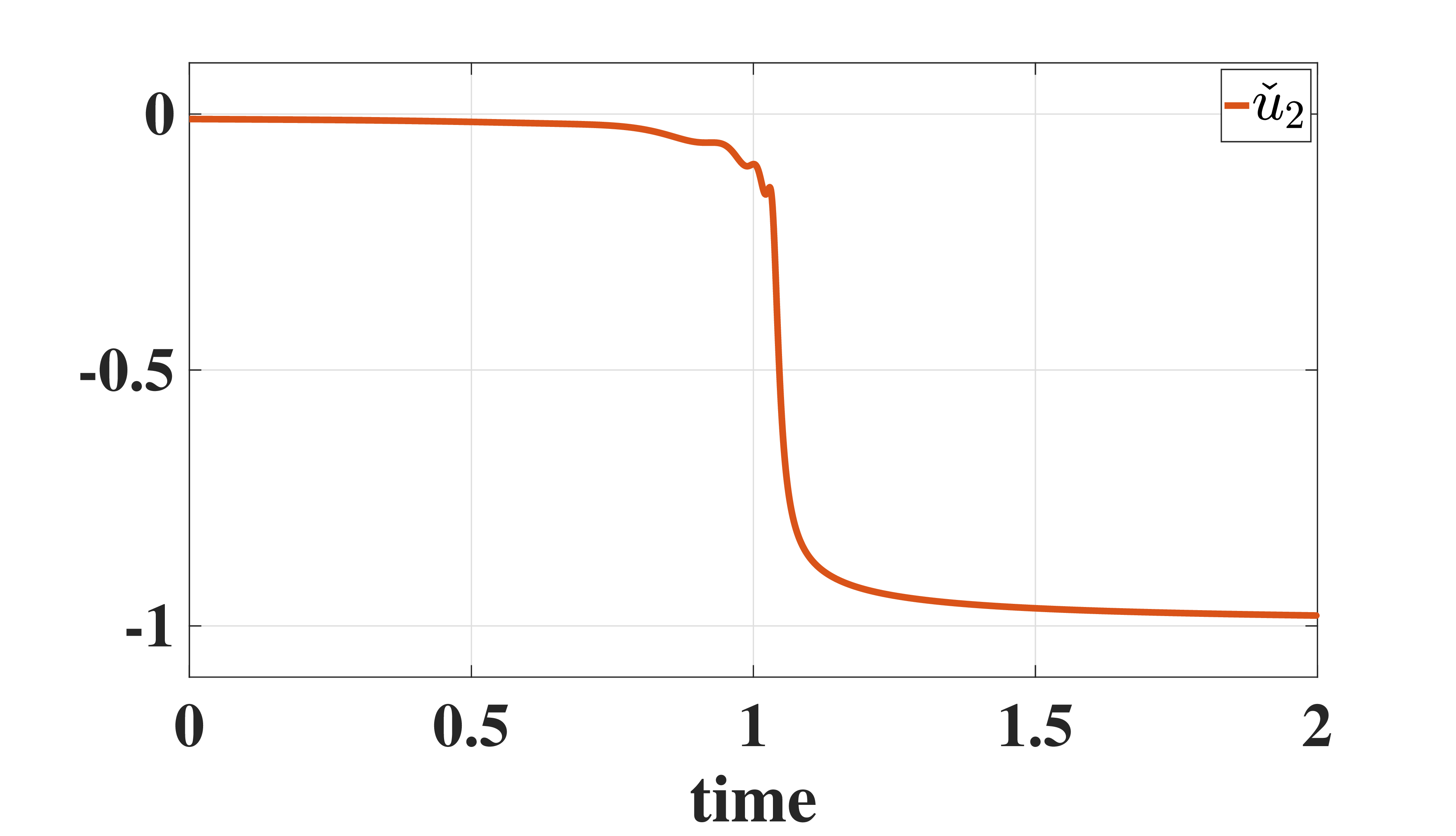}
	\includegraphics[width=0.49\linewidth]{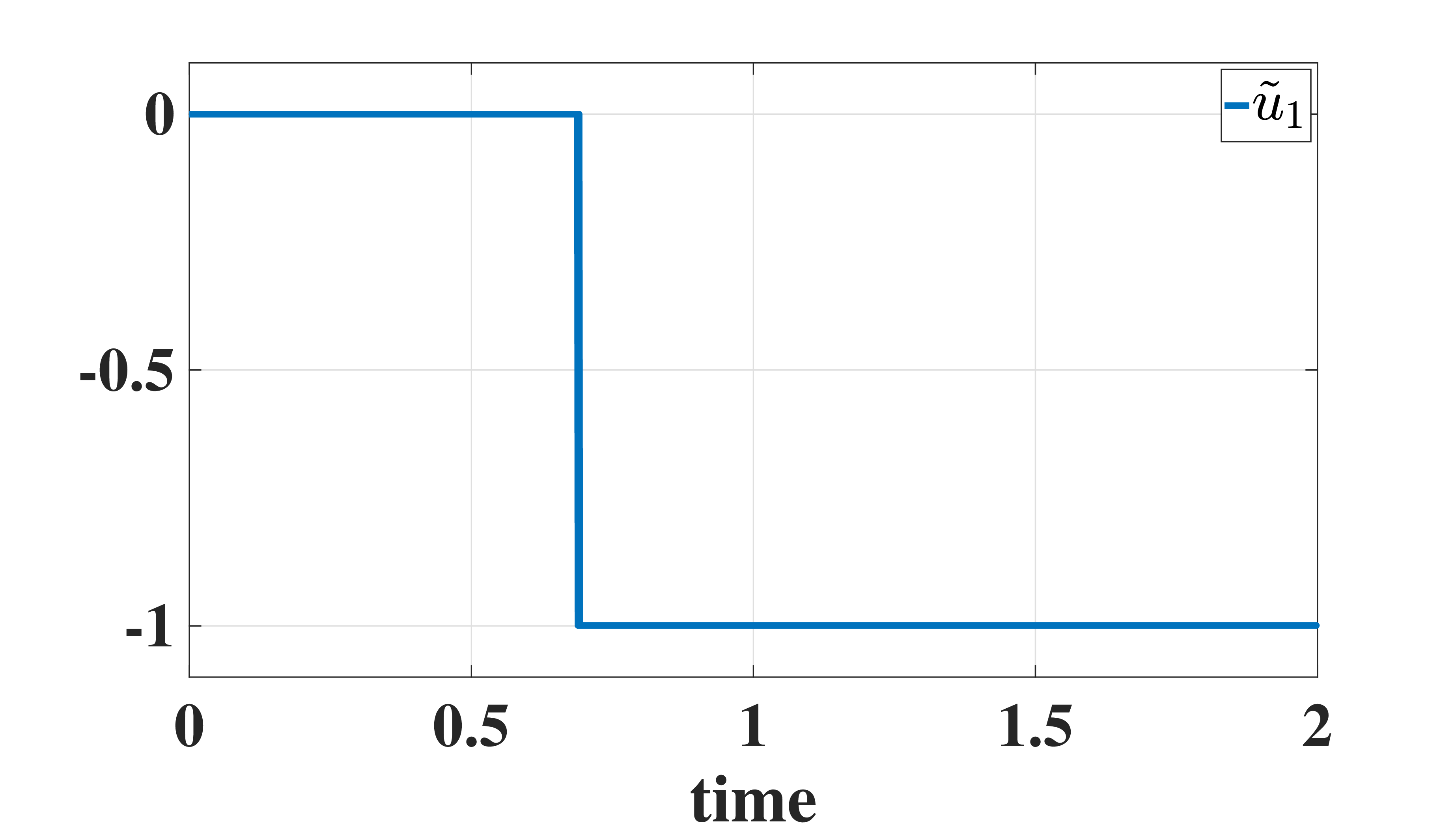}
	\includegraphics[width=0.49\linewidth]{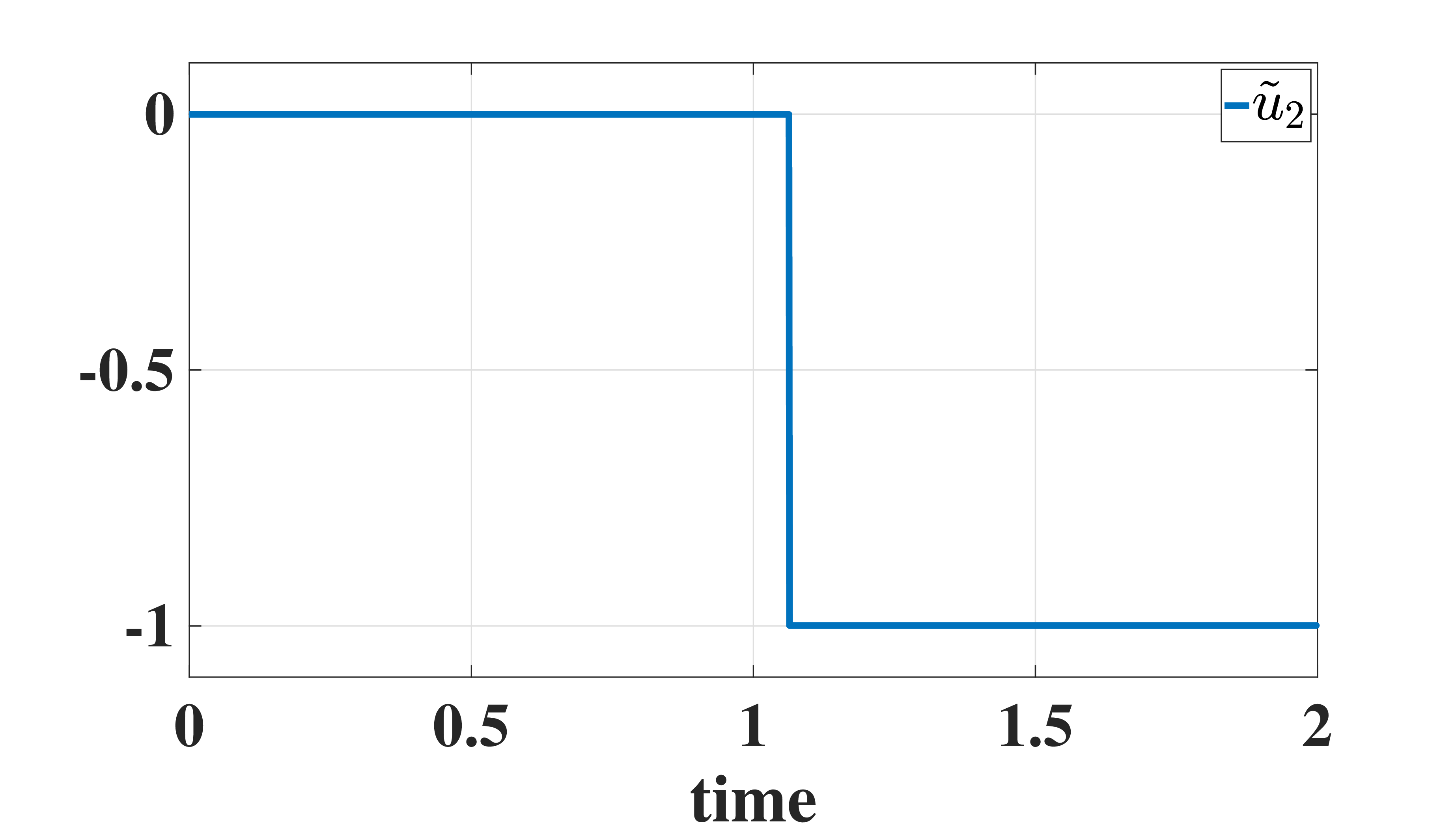}
		\caption{The $L^1$ optimal control $\check{u}$ (top) and the non-convex optimal control $\tilde{u}$ (bottom) in Example~2.}
	\label{fig:ex2_nonconvex}
\end{figure}

For the sparse optimal control problem with the parameters
$x_0(\alpha) = \sin(2\pi\alpha)$, $x_f = 0$, $T = 2$, and $\lambda = 10^6$,
Figure~\ref{fig:ex2_nonconvex} shows the $L^1$ optimal control $\check{u}$ and the non-convex optimal control $\tilde{u}$, where the non-convex penalty function is taken as
\begin{align*}
	&\psi(u) =\frac{10}{9}\left(\left\|u\right\|_{\ell^1} - 0.1 \left\|u\right\|_{\ell^2}^2\right).
\end{align*}
In both cases, the terminal state error was $0.02194$, confirming that the desired final state was nearly achieved.  
Furthermore, this slight modification of the cost function, where a small multiple of the squared $\ell^2$ norm is subtracted from the $\ell^1$ norm, is shown to effectively induce the sparsity in the resulting control.

\subsection{Example 3}\label{subsec:example_3}

We consider the finite-dimensional system~\eqref{eq:system_finite} where 
the underlying graph $G^{[n]}$ is given by a simple threshold graph\footnote{The simple threshold graph is defined on the set $\{1,2,  \dots, n\}$ by connecting $i$ and $j$ if and only if $i + j \leq n+1$.}, 
the matrix $\mathcal{B}^{[n]}$ is given by
\begin{align*}
	\mathcal{B}_1^{[n]} = \sum_{i=1}^{i_1} e_i^{[n]},\quad
	\mathcal{B}_2^{[n]} = \sum_{i= i_1 + 1}^{i_2} e_i^{[n]},\quad
	\mathcal{B}_3^{[n]} = \sum_{i= i_2+1}^{n} e_i^{[n]},
\end{align*}
where $i_1 = \lceil \frac{n}{4}\rceil$, $i_2 = \lceil \frac{3n}{4}\rceil$ ($\lceil \cdot \rceil$ is the ceiling function), and $e_i^{[n]}$ is the $i$th canonical vector in $\mathbb{R}^{n}$.
The initial state is given by $\mathrm{x}^{[n]}_0=0$
and the target state $\mathrm{x}_f^{[n]}$ is given by 
$\mathrm{x}_{f,i}^{[n]} = n \int_{P_i} \sqrt{1-\alpha^2} d\alpha$
for $i\in\left\llbracket n \right\rrbracket$.
Note that this system satisfies Assumption~\ref{ass:convergence} for 
\begin{align*}
	&W(\alpha, \beta) = {\bf 1}_{\Omega}(\alpha, \beta), \\
	&\Omega = \{(\alpha, \beta) \in [0,1]^2: 0 \leq \alpha + \beta \leq 1\},\\
	&b_1(\alpha) = {\bf 1}_{\Omega_1}(\alpha),\quad
		b_2(\alpha) = {\bf 1}_{\Omega_2}(\alpha),\quad
		b_3(\alpha) = {\bf 1}_{\Omega_3}(\alpha),\\
	&\Omega_1 = [0,0.25), \quad
		\Omega_2 = [0.25, 0.75),\quad
		\Omega_3 = [0.75, 1],\\
	&x_0 = 0,\quad
		x_f(\alpha) = \sqrt{1-\alpha^2}.
\end{align*}
Note also that the operator $A = T_W$ has eigenvalues $\lambda_k = \frac{2}{(4k+1)\pi}$ and eigenvectors $\sqrt{2}\cos\frac{\alpha}{\lambda_k}$, where $k\in\mathbb{Z}$; the computation follows the method in~\cite[Example A.4.20]{CurZwa20}.
Hence, this example satisfies Assumption~\ref{ass:BOB}.

\begin{figure}[t]
	\centering
	\includegraphics[width=0.8\linewidth]{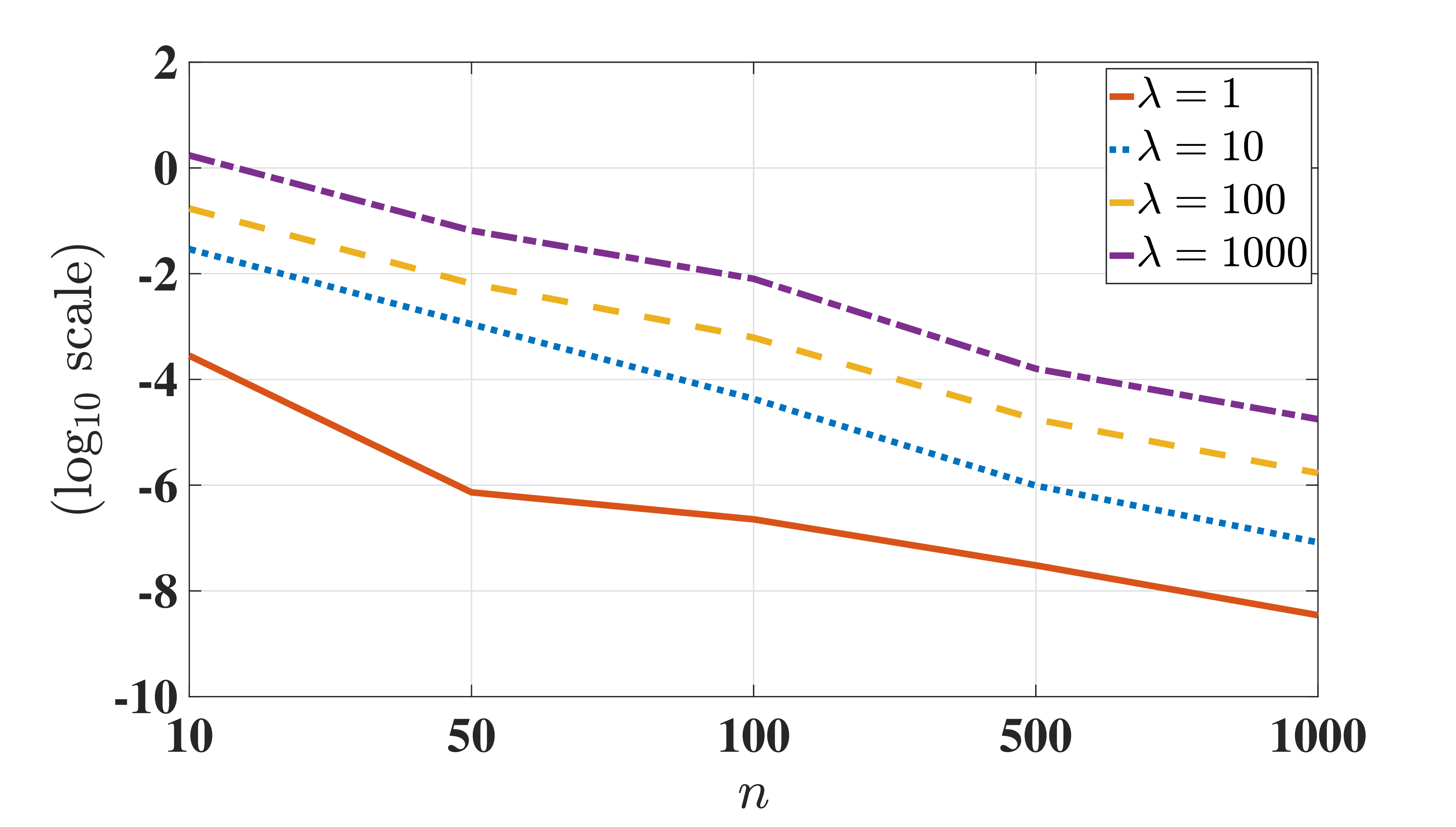}
	\caption{Approximation error $J_0^{[n]}(\bar{u}) - \bar{J}_0^{[n]}$ in Example~3.}
	\label{fig:error_ex3}
\end{figure}

We consider Problem~\ref{prob:sparse_finite} for the parameters
$T = 1$ and $\lambda\in\{10^0, 10^1, 10^2, 10^3\}$.
Based on Theorem~\ref{thm:L0-L1}, we first computed the sparse optimal control $\bar{u}$ by solving Problem~\ref{prob:convex_limit} for the limit graphon system $(A; B; x_0)$.
Figure~\ref{fig:error_ex3} shows the approximation errors $J_0^{[n]}(\bar{u}) - \bar{J}_0^{[n]}$ for $n\in\{10, 50, 100, 500, 1000\}$, 
where $\bar{J}_0^{[n]}$ is the optimal value of Problem~\ref{prob:sparse_finite} and this is equal to the optimal value of the corresponding $L^1$ optimal control problem~\cite{IkeKas_ECC19}.
The errors asymptotically approach zero as the number of nodes $n$ increases, as established in Theorem~\ref{thm:app_convergence}.

\section{Conclusion}\label{sec:conclusion}

This paper has addressed a sparse optimal control problem for infinite-dimensional systems.
We have shown two results: 
(i) theoretical conditions under which the sparse optimal control can be obtained from the $L^1$ optimal control problem and a class of non-convex optimal control problems, 
and (ii) an approximation method for sparse optimal controls in finite-dimensional systems based on the graph limit theory.
This paper has focused on a feedforward optimal control for a given initial state. 
On the other hand, it is also important to consider the influences of disturbances and variations in the size of the finite-dimensional system under consideration.
Future work therefore includes extending the proposed approach to a state-feedback control method, such as model predictive control.
To this end, the authors are currently pursuing approximate state-feedback control laws based on a dynamic programming approach.

\bibliographystyle{plain}        
\bibliography{autosam}           

\appendices

\section{Proof of Proposition~\ref{prop:ass_condition}}
\label{sec:app_assumption}

We first assume that the system $(A; \mathcal{M})$ is approximately controllable on $[0, T]$.
From~\cite[Theorem 6.2.6]{CurZwa20}, 
\begin{equation}\label{eq:app_contro_1}
	\mathcal{M}^{\prime} {e^{At}}^\prime z = 0
\end{equation}
on $[0, T]$ implies $z=0$.
Since $A$ is self-adjoint, we have
\begin{align*}
	\mathcal{M}^{\prime} {e^{At}}^\prime z (1, 1) 
	=\left(\mathcal{M}^{\prime} {e^{At}}^\prime z, 1\right) 
	= \left({e^{At}}^\prime z, \mathcal{M} 1\right) 
	= \left({e^{At}}^\prime z, Ab\right)
	= \left(z, {e^{At}} Ab\right) 
	= \left(z, A {e^{At}} b\right) 
	= \left(Az, {e^{At}} b\right).
\end{align*}
Hence, \eqref{eq:app_contro_1} holds if and only if $\left(Az, {e^{At}} b\right) = 0$.
This implies that $0$ is not an eigenvalue of $A$. 
(Indeed, if $0$ were an eigenvalue of $A$, then its corresponding eigenvector $z \neq 0$ satisfies \eqref{eq:app_contro_1}.)

Since $A\in\mathcal{L}(L^2_{[0,1]})$ is compact and self-adjoint, there exists an orthonormal basis of $L^2_{[0,1]}$ consisting of eigenvectors $\{\phi_{i_k}\}$ such that
\[
	A z = \sum_{i=1}^{\infty} \rho_i \sum_{k=1}^{r_i} \left(z, \phi_{i_k}\right) \phi_{i_k}
\]
holds for all $z \in L^2_{[0,1]}$~\cite[Theorem A.4.19]{CurZwa20},
where
$\{\rho_i\} \subset [-\|A\|_{\mathrm{op}}, \|A\|_{\mathrm{op}}]$ are distinct eigenvalues of $A$,
$\{\phi_{i_k}\}$, $k = 1, 2, \dots, r_i$, are eigenvectors corresponding to $\rho_i$,
and $r_i \geq 1$ is the (finite) dimension of the eigenspace corresponding to $\rho_i$.
From~\cite[Theorem 6.2.25 and Theorem 6.3.4]{CurZwa20}, it follows that
\begin{equation}\label{eq:app_contro_2}
	\mathrm{rank~} 
	\begin{bmatrix} \left(Ab, \phi_{i_1}\right) & \left(Ab, \phi_{i_2}\right) & \cdots & \left(Ab, \phi_{i_{r_i}}\right) \end{bmatrix} = r_i
\end{equation}
for all $i$.
Since we have
\begin{align*}
	\left(Ab, \phi_{i_k}\right) 
	= \left(b, A\phi_{i_k}\right)
	=\rho_i \left(b, \phi_{i_k}\right)
\end{align*} 
and $\rho_i \neq 0$,
\eqref{eq:app_contro_2} implies $r_i=1$ and $\left(b, \phi_{i}\right)\neq0$ for all $i$, where $\phi_i$ is a normalized eigenvector corresponding to $\rho_i$.

We next assume that 
$\lambda_i \neq 0$ for all $i$, 
$\lambda_i \neq \lambda_j$ for all $i\neq j$, 
and $\left(b, \phi_i\right) \neq 0$ for all $i$.
Then, \eqref{eq:app_contro_2} holds, and the system $(A; \mathcal{M})$ is approximately controllable on $[0, T]$ from~\cite[Theorem 6.2.25 and Theorem 6.3.4]{CurZwa20}.

\section{Proof of Theorem~\ref{thm:discrete-convex}}\label{sec:app_analytic}

\begin{lemma}\label{lem:analytic}
For $A\in\mathcal{L}(L^2_{[0,1]})$ and $a, b \in L^2_{[0,1]}$, 
define a function $\omega: [0, T] \to \mathbb{R}$ by
$\omega(t) = \left( a, e^{A(T-t)} b \right)$.
Then, $\omega$ is continuous on $[0,T]$ and analytic on $(0, T)$, 
and we have
\begin{equation}\label{eq:theta_dif}
	\frac{d^k \omega}{dt^k}(t) = (-1)^k \left( a, e^{A(T-t)} A^k b \right)
\end{equation}
for any $k\in\mathbb{N}$.
\end{lemma}

\begin{IEEEproof}
The continuity of $\omega$ immediately follows from the continuity of the functions $(a, \cdot)$ and $e^{At} b$.
%
From~\cite[Theorem 2.1.13]{CurZwa20}, \eqref{eq:theta_dif} holds for $k=1$.
For $k\geq 2$, \eqref{eq:theta_dif} can be shown by induction.
For any $t\in(0, T)$, we have
\begin{align*}
	\left|\frac{d^k \omega}{dt^k}(t)\right|
	\leq \left\| a \right\|_2 \left\|e^{A(T-t)} A^k b\right\|_2
	\leq \left\| a \right\|_2 e^{\|A\|_{\mathrm{op}}T} \|A\|_{\mathrm{op}}^{k} \left\|b\right\|_2
\end{align*}
for any $k\in\mathbb{N}$.
From~\cite[Proposition 1.2.12]{Kra02}, $\omega$ is analytic on $(0, T)$.
\end{IEEEproof}

\section{Proof of Lemma~\ref{lem:app_limit}}\label{sec:app_convergence}

\begin{lemma}\label{lem:convergence}
Suppose that Assumption~\ref{ass:convergence} holds. 
Define $\pi_1^{[n]}$, $\pi_2^{[n]}$, $\pi_3$, $\pi_4$, $\pi_5^{[n]}$, $\pi_6^{[n]}\in L^2_{[0,1]}$ by \eqref{def:pi} and \eqref{def:pi_statement3}.
Then, the following hold:
\begin{align}
	&\lim_{n\to\infty} \left\| \pi_1^{[n]} - \pi_3 \right\|_{2} = 0,\label{lem_eq:lemma_1}\\
	&\lim_{n\to\infty} \left\| \pi_2^{[n]} -\pi_4 \right\|_{2} = 0, \label{lem_eq:lemma_2}\\
	&\lim_{n\to\infty} \left(\pi_1^{[n]}, \pi_2^{[n]}\right) = \left(\pi_3, \pi_4\right),\label{lem_eq:lemma_3}\\
	&\lim_{n\to\infty} \left(\pi_1^{[n]}, x_f^{[n]}\right) = \left(\pi_3, x_f\right),\label{lem_eq:lemma_3_2}\\
	&\lim_{n\to\infty} \left(\pi_2^{[n]}, x_f^{[n]}\right) = \left(\pi_4, x_f\right),\label{lem_eq:lemma_3_3}\\
	&\lim_{n\to\infty} \left\| \pi_5^{[n]} -\pi_6^{[n]} \right\|_{2} = 0, \label{lem_eq:lemma_4}\\
	&\lim_{n\to\infty} \left| \left(\pi_1^{[n]}, \pi_5^{[n]}\right) - \left(\pi_3, \pi_6^{[n]}\right) \right| = 0, \label{lem_eq:lemma_5}\\
	&\lim_{n\to\infty} \left| \left(\pi_6^{[n]}, x_f \right) - \left(\pi_5^{[n]}, x_f^{[n]} \right)\right| = 0.\label{lem_eq:lemma_6}
\end{align}
\end{lemma}

\begin{IEEEproof}
Since we have
\begin{align*}
	\left\|\pi_1^{[n]} - \pi_3\right\|_{2}
	&\leq \left\|e^{A^{[n]}T}x_{0}^{[n]} - e^{A^{[n]}T}x_{0} \right\|_{2}
		+ \left\|e^{A^{[n]}T}x_{0} - e^{AT}x_{0} \right\|_{2}\\
	&\leq e^{\left\|A^{[n]} - A\right\|_{\mathrm{op}}T}  e^{\left\|A\right\|_{\mathrm{op}}T}  \left\|x_{0}^{[n]} - x_{0} \right\|_{2}
		+ e^{\left\|A\right\|_{\mathrm{op}}T} \left(e^{\left\|A^{[n]} - A\right\|_{\mathrm{op}}T} - 1\right) \left\|x_{0} \right\|_{2},
\end{align*}
we obtain \eqref{lem_eq:lemma_1}, where the second inequality follows from \cite[p. 78]{Paz12}.
Since we have
\begin{align*}
	\left\| \pi_2^{[n]} - \pi_4 \right\|
	&\leq \int_{0}^{T} \left\| e^{A^{[n]}(T-t)} B^{[n]} \bar{u}(t) - e^{A(T-t)} B \bar{u}(t) \right\|_{2} dt\\
	&\leq \int_{0}^{T} \left\| e^{A^{[n]}(T-t)} B^{[n]} \bar{u}(t) - e^{A^{[n]}(T-t)} B \bar{u}(t) \right\|_{2} dt 
		+ \int_{0}^{T} \left\| e^{A^{[n]}(T-t)} B \bar{u}(t) - e^{A(T-t)} B \bar{u}(t)\right\|_{2} dt \\
	&\leq \int_{0}^{T} e^{\left\|A^{[n]}\right\|_{\mathrm{op}}(T-t)} \left\|B^{[n]} \bar{u}(t) - B \bar{u}(t)\right\|_{2} dt
		+ \int_{0}^{T} e^{\left\|A\right\|_{\mathrm{op}}(T-t)} \left(e^{\left\|A^{[n]} - A\right\|_{\mathrm{op}}(T-t)} - 1\right) \left\|B \bar{u}(t) \right\|_{2} dt \\
	&\leq T e^{\left\|A^{[n]} - A\right\|_{\mathrm{op}}T} e^{\left\|A\right\|_{\mathrm{op}}T} \sum_{j=1}^{m} \left\|b_j^{[n]} - b_j \right\|_{2}
		+ T e^{\left\|A\right\|_{\mathrm{op}}T} \left(e^{\left\|A^{[n]} - A\right\|_{\mathrm{op}}T} - 1 \right) \sum_{j=1}^{m} \left\| b_j \right\|_{2},		
\end{align*}
we obtain \eqref{lem_eq:lemma_2}. 
Since we have
\begin{align*}
	\left(\pi_1^{[n]}, \pi_2^{[n]}\right) - \left(\pi_3, \pi_4\right)
	= \left(\pi_1^{[n]} - \pi_3, \pi_4\right) + \left(\pi_1^{[n]} - \pi_3, \pi_2^{[n]} - \pi_4\right) + \left(\pi_3, \pi_2^{[n]} - \pi_4\right),
\end{align*}
we obtain \eqref{lem_eq:lemma_3} from Cauchy-Schwarz inequality, \eqref{lem_eq:lemma_1}, and \eqref{lem_eq:lemma_2}.
Also, we obtain \eqref{lem_eq:lemma_3_2} and \eqref{lem_eq:lemma_3_3} in a similar evaluation. 
Since we have
\begin{align*}
	\left\| \pi_5^{[n]} - \pi_6^{[n]} \right\|_{2}
	&\leq \int_{0}^{T} \left\| e^{A^{[n]}(T-t)} B^{[n]} \bar{u}^{[n]}(t) - e^{A(T-t)} B \bar{u}^{[n]}(t) \right\|_{2} dt\\
	&\leq \int_{0}^{T} \left\| e^{A^{[n]}(T-t)} B^{[n]} \bar{u}^{[n]}(t) - e^{A^{[n]}(T-t)} B \bar{u}^{[n]}(t) \right\|_{2} dt
		+ \int_{0}^{T} \left\| e^{A^{[n]}(T-t)} B \bar{u}^{[n]}(t) - e^{A(T-t)} B \bar{u}^{[n]}(t) \right\|_{2} dt \\
	&\leq T e^{\left\|A^{[n]} - A\right\|_{\mathrm{op}}T} e^{\left\|A\right\|_{\mathrm{op}}T} \sum_{j=1}^{m} \left\|b_j^{[n]} - b_j \right\|_{2}
		+ T e^{\left\|A\right\|_{\mathrm{op}}T} \left(e^{\left\|A^{[n]} - A\right\|_{\mathrm{op}}T} - 1\right) \sum_{j=1}^{m} \left\|b_j\right\|_{2},
\end{align*}
we obtain \eqref{lem_eq:lemma_4}. 
Since we have
\begin{align*}
	\left| \left(\pi_1^{[n]}, \pi_5^{[n]}\right) - \left(\pi_3, \pi_6^{[n]}\right) \right|
	\leq \left\|\pi_1^{[n]} - \pi_3\right\|_2 \left\|\pi_5^{[n]}\right\|_2 + \left\|\pi_3\right\|_2 \left\|\pi_5^{[n]} - \pi_6^{[n]} \right\|_2
\end{align*}
and
\begin{align*}
	\left\| \pi_5^{[n]} \right\|_{2}
	\leq T e^{\left\|A^{[n]} - A \right\|_{\mathrm{op}}T} e^{\left\|A\right\|_{\mathrm{op}}T}  \sum_{j=1}^{m} \left(\left\|b_{j}^{[n]}-b_j\right\|_{2} + \left\|b_{j}\right\|_{2}\right),
\end{align*}
we obtain~\eqref{lem_eq:lemma_5} from~\eqref{lem_eq:lemma_1} and \eqref{lem_eq:lemma_4}.
Also, we obtain \eqref{lem_eq:lemma_6} in a similar evaluation.
\end{IEEEproof}


\end{document}